\documentclass[11pt]{article}

\usepackage{amssymb}
\usepackage{amsfonts}
\usepackage{amsmath}
\usepackage{amsthm}
\usepackage[margin=0.85in]{geometry}
\usepackage{enumerate}
\usepackage{amstext}
\usepackage{layout}

\numberwithin{equation}{section}

\newtheorem{theorem}{Theorem}[section]
\newtheorem{corollary}{Corollary}[section]
\newtheorem{lemma}{Lemma}[section]

\newtheorem{remark}{Remark}[section]
\newtheorem{definition}{Definition}[section]
\newtheorem{example}{Example}[section]

\begin{document}

\noindent {\bf\large{A Characterization of a New Type of Strong Law of
Large Numbers}}

\vskip 0.3cm

\noindent {\bf Deli Li\footnote{Deli Li, Department of Mathematical
Sciences, Lakehead University, Thunder Bay, Ontario, Canada P7B
5E1\\ e-mail: dli@lakeheadu.ca} $\cdot$ Yongcheng
Qi\footnote{Yongcheng Qi, Department of Mathematics and
Statistics, University of Minnesota Duluth, Duluth, Minnesota 55812, U.S.A.\\
e-mail: yqi@d.umn.edu}$\cdot$ Andrew Rosalsky\footnote{Andrew
Rosalsky, Department of Statistics, University of Florida,
Gainesville, Florida 32611, U.S.A.\\ e-mail: rosalsky@stat.ufl.edu}
\footnote{Corresponding author: Andrew Rosalsky (Telephone:
1-352-273-2983, FAX: 1-352-392-5175)}}

\vskip 0.3cm

\noindent {\bf Abstract}~~Let $0 < p < 2$ and $1 \leq q < \infty$.
Let $ \{X_{n};~n \geq 1 \}$ be a sequence of independent copies
of a real-valued random variable $X$ and set
$S_{n} = X_{1} + \cdots + X_{n}, ~n \geq 1$.
We say $X$ satisfies the $(p, q)$-{\it type strong law of large numbers}
(and write $X \in SLLN(p, q)$) if
$\sum_{n = 1}^{\infty} \frac{1}{n}\left(\frac{\left|S_{n}\right|}{n^{1/p}}
\right)^{q} < \infty$ almost surely. This paper is devoted to a characterization
of $X \in SLLN(p, q)$. By applying results obtained from the new versions of 
the classical L\'{e}vy, Ottaviani, and Hoffmann-J{\o}rgensen (1974) inequalities 
proved by Li and Rosalsky (2013) and by using techniques developed by Hechner 
and Heinkel (2010), we show that $X \in SLLN(p, q)$ if and only if
\[
\left \{
\begin{array}{ll}
\mbox{$\displaystyle \mathbb{E}X = 0 ~~\mbox{and}~~\int_{0}^{\infty}
\mathbb{P}^{q/p}\left(|X|^{q}
> t \right) dt < \infty$} & \mbox{if $1 \leq q < p < 2$,}\\
&\\
\mbox{$\displaystyle \mathbb{E}X = 0, ~~\mathbb{E} |X|^{p} < \infty, ~~\mbox{and}
~~\sum_{n=1}^{\infty} \frac{\int_{\min\{u_{n}^{p},
n\}}^{n} \mathbb{P}\left(|X|^{p} > t \right) dt}{n} < \infty$} & \mbox{if
~~$\displaystyle 1 < p = q < 2$,}\\
&\\
\mbox{$\displaystyle \mathbb{E}X = 0 ~~\mbox{and}~~\mathbb{E}|X|^{p} < \infty$}
& \mbox{if $1 < p < 2$ \mbox{and} $q > p$,}\\
&\\
\mbox{$\displaystyle \mathbb{E}X = 0, ~~\sum_{n=1}^{\infty}
\frac{|\mathbb{E}XI\left\{|X| \leq n \right\}|}{n} < \infty,~\mbox{and}$} &\\
&\\
\mbox{$\displaystyle \sum_{n=1}^{\infty} \frac{\int_{\min\{u_{n}, n\}}^{n}
\mathbb{P}\left(|X| > t \right) dt}{n} < \infty$} & \mbox{if
~~$\displaystyle q = p = 1$,}\\
&\\
\mbox{$\displaystyle \mathbb{E}X = 0~~\mbox{and}~~\sum_{n=1}^{\infty}
\frac{\left|\mathbb{E}XI\left\{|X| \leq n \right\} \right|^{q}}{n} < \infty$}
& \mbox{if $p = 1 < q$,}\\
&\\
\mbox{$\displaystyle \mathbb{E}|X|^{p} < \infty$}
& \mbox{if $0 < p < 1 \leq q$,}
\end{array}
\right.
\]
where $u_{n} = \inf \left\{t:~ \mathbb{P}(|X| > t) < \frac{1}{n}
\right\}, ~n \geq 1$. For $q = 1$, this equivalence has recently been discovered by Li, Qi, and
Rosalsky (2011). Versions of above results in a Banach space setting are
also presented.

~\\


\vskip 0.3cm


\vskip 0.3cm

\noindent {\bf Running Head}: Strong law of large numbers

\section{Introduction}

Throughout, let $(\mathbf{B}, \| \cdot \| )$ be a real separable
Banach space equipped with its Borel $\sigma$-algebra $\mathcal{B}$
($=$ the $\sigma$-algebra generated by the class of open subsets of
$\mathbf{B}$ determined by $\|\cdot\|$) and let $ \{X_{n};~n \geq 1
\}$ be a sequence of independent copies of a {\bf B}-valued random
variable $X$ defined on a probability space $(\Omega, \mathcal{F},
\mathbb{P})$. As usual, let $S_{n} = \sum_{k=1}^{n} X_{k},~ n \geq
1$ denote their partial sums. If $0 < p < 2$ and if $X$ is a
real-valued random variable (that is, if $\mathbf{B} = \mathbb{R}$),
then
\[
\lim_{n \rightarrow \infty} \frac{S_{n}}{n^{1/p}} = 0 ~~\mbox{almost
surely (a.s.)}
\]
if and only if
\[
\mathbb{E}|X|^{p} < \infty ~~\mbox{where}~~ \mathbb{E}X = 0
~~\mbox{whenever}~~ p \geq 1.
\]
This is the celebrated Kolmogorov-Marcinkiewicz-Zygmund strong law
of large numbers (SLLN); see Kolmogoroff [8] for $p = 1$ and
Marcinkiewicz and Zygmund [13] for $p \neq 1$.

The classical Kolmogorov SLLN in real separable Banach spaces was
established by Mourier [17]. The extension of the
Kolmogorov-Marcinkiewicz-Zygmund SLLN to $\mathbf{B}$-valued random
variables is independently due to Azlarov and Volodin [1] and de
Acosta [3].

\vskip 0.2cm

\begin{theorem}
{\rm (Azlarov and Volodin [1] and de Acosta [3])}. Let $0 < p < 2$
and let $\{X_{n}; ~n \geq 1\}$ be a sequence of independent copies
of a $\mathbf{B}$-valued random variable $X$. Then
\[
\lim_{n \rightarrow \infty} \frac{S_{n}}{n^{1/p}} = 0~~\mbox{a.s.}
\]
if and only if
\[
\mathbb{E}\|X\|^{p} < \infty~~\mbox{and}~~\frac{S_{n}}{n^{1/p}}
\rightarrow_{\mathbb{P}} 0.
\]
\end{theorem}

\vskip 0.2cm

Let $\{R_{n}; ~n \geq 1\}$ be a {\it Rademacher sequence}; that is,
$\{R_{n};~n \geq 1\}$ is a sequence of independent and identically
distributed (i.i.d.) random variables with $\mathbb{P}\left(R_{1} =
1\right) = \mathbb{P}\left(R_{1} = -1\right) = 1/2$. Let
$\mathbf{B}^{\infty} =
\mathbf{B}\times\mathbf{B}\times\mathbf{B}\times \cdots$ and define
\[
\mathcal{C}(\mathbf{B}) = \left\{(v_{1}, v_{2}, ...) \in
\mathbf{B}^{\infty}: ~\sum_{n=1}^{\infty} R_{n}v_{n}
~~\mbox{converges in probability} \right\}.
\]
Let $1 \leq p \leq 2$. Then $\mathbf{B}$ is said to be of {\it
Rademacher type $p$} if there exists a constant $0 < C < \infty$
such that
\[
\mathbb{E}\left\|\sum_{n=1}^{\infty}R_{n}v_{n} \right\|^{p} \leq C
\sum_{n=1}^{\infty}\|v_{n}\|^{p}~~\mbox{for all}~(v_{1}, v_{2}, ...)
\in \mathcal{C}(\mathbf{B}).
\]
Hoffmann-J{\o}rgensen and Pisier [7] proved for $1 \leq p \leq 2$
that $\mathbf{B}$ is of Rademacher type $p$ if and only if there
exists a constant $0 < C < \infty$ such that
\[
\mathbb{E}\left\|\sum_{k=1}^{n}V_{k} \right\|^{p} \leq C
\sum_{k=1}^{n} \mathbb{E}\left\|V_{k}\right\|^{p}
\]
for every finite collection $\{V_{1}, ..., V_{n} \}$ of independent
mean $0$ $\mathbf{B}$-valued random variables.

If $\mathbf{B}$ is of Rademacher type $p$ for some $p \in (1, 2]$,
then it is of Rademacher type $q$ for all $q \in [1, p)$. Every real
separable Banach spaces is of Rademacher type (at least) $1$.

Let $0 < p \leq 2$ and let $\{\Theta_{n}; ~n \geq 1 \}$ be a
sequence of i.i.d. stable random variables each with characteristic
function $\psi(t) = \exp \left\{-|t|^{p}\right\}, ~- \infty < t <
\infty$. Then $\mathbf{B}$ is said to be of {\it stable type $p$} if
$\sum_{n=1}^{\infty} \Theta_{n}v_{n}$ converges a.s. whenever
$\{v_{n}: ~n \geq 1\} \subseteq \mathbf{B}$ with
$\sum_{n=1}^{\infty} \|v_{n}\|^{p} < \infty$. Equivalent
characterizations of a Banach space being of stable type $p$,
properties of stable type $p$ Banach spaces, as well as various
relationships between the conditions ``Rademacher type $p$" and
``stable type $p$" may be found in Maurey and Pisier [16],
Woyczy\'{n}ski [20], Marcus and Woyczy\'{n}ski [15], Rosi\'{n}ski
[19], Pisier [18], and Ledoux and Talagrand [9]. Some of these properties
and relationships are summarized in Li, Qi, and Rosalsky [10].

De Acosta [3] also provided a remarkable characterization of
Rademacher type $p$ Banach spaces. Specifically, de Acosta [3]
proved the following theorem.

\vskip 0.2cm

\begin{theorem}
{\rm (de Acosta [3])}. Let $1 \leq p < 2$. Then the following two
statements are equivalent:
\begin{align*}
& {\bf (i)} \quad \mbox{The Banach space $\mathbf{B}$ is of
Rademacher type $p$.}\\
& {\bf (ii)} \quad \mbox{For every sequence $\{X_{n}; ~n \geq 1 \}$
of independent copies of a $\mathbf{B}$-valued variable $X$},
\end{align*}
\[
\lim_{n \rightarrow \infty} \frac{S_{n}}{n^{1/p}} = 0~~\mbox{a.s. if
and only if}~~\mathbb{E}\|X\|^{p} < \infty~~\mbox{and}~~\mathbb{E}X
= 0.
\]
\end{theorem}

\vskip 0.2cm

At the origin of the current investigation are the following recent
and striking result by Hechner and Heinkel [5] which is new even in the
case where the Banach space $\mathbf{B}$ is the real line. The earliest
investigation that we are aware of concerning the convergence of the
series $\sum_{n=1}^{\infty} \frac{1}{n} \left(\frac{\mathbb{E}|S_{n}|}{n} \right)$
was carried out by Hechner [4] for the case where $\{X_{n};~n \geq 1\}$
is a sequence of i.i.d. mean $0$ real-valued random variables.

\vskip 0.2cm

\begin{theorem}
{\rm (Hechner and Heinkel [5])}. Suppose that $\mathbf{B}$ is of
stable type $p$ ($1 < p < 2$) and let $\{X_{n}; ~n \geq 1 \}$ be a
sequence of independent copies of a $\mathbf{B}$-valued variable $X$
with $\mathbb{E}X = 0$. Then
\[
\sum_{n=1}^{\infty} \frac{1}{n}
\left(\frac{\mathbb{E}\|S_{n}\|}{n^{1/p}} \right) < \infty
\]
if and only if
\[
\int_{0}^{\infty} \mathbb{P}^{1/p}(\|X\| > t) dt < \infty.
\]
\end{theorem}

\vskip 0.2cm

Inspired by the above discovery by Hechner and Heinkel [5], Li, Qi,
and Rosalsky [10] obtained sets of necessary and sufficient conditions
for
\[
\sum_{n=1}^{\infty} \frac{1}{n} \left(\frac{\|S_{n}\|}{n^{1/p}}
\right) < \infty~~\mbox{a.s.}
\]
for the three cases: $0 < p < 1$, $p = 1$, $1 < p < 2$. Moreover, Li, Qi,
and Rosalsky [10] obtained necessary and sufficient conditions for
\[
\sum_{n=1}^{\infty} \frac{1}{n} \left(\frac{\mathbb{E}\|S_{n}\|}{n}
\right) < \infty.
\]
Again, these results are new when $\mathbf{B} = \mathbb{R}$; see
Theorem 2.5 of Li, Qi, and Rosalsky [10].

Motivated by the results obtained by Li, Qi, and Rosalsky [10], we introduce
a new type strong law of large numbers as follows.

\vskip 0.2cm

\begin{definition}
Let $0 < p < 2$ and $0 < q < \infty$. Let $ \{X_{n};~n \geq 1 \}$ be a sequence of
independent copies of a $\mathbf{B}$-valued random variable $X$. We say $X$ satisfies
the $(p, q)$-type strong law of large numbers (and write $X \in SLLN(p, q)$) if
\[
\sum_{n = 1}^{\infty} \frac{1}{n}\left(\frac{\left\|S_{n}\right\|}{n^{1/p}}
\right)^{q} < \infty ~~\mbox{a.s.}
\]
\end{definition}

The following result was recently obtained by Li, Qi, and Rosalsky [11] who
proved it by employing new versions of the classical L\'{e}vy, Ottaviani,
and Hoffmann-J{\o}rgensen [6] inequalities established by Li and Rosalsky [12]
and by using some of techniques developed by Hechner and Heinkel [5]. Note
that no conditions are imposed on the Banach space $\mathbf{B}$. Theorem 1.4
will be used in the proofs of the main results of the current work.

\vskip 0.3cm

\begin{theorem}
Let $0 < p < 2$ and $0 < q < \infty$. Let $ \{X_{n};~n \geq 1 \}$ be a sequence
of independent copies of a $\mathbf{B}$-valued random variable $X$. Then
\begin{equation}
\sum_{n = 1}^{\infty} \frac{1}{n} \mathbb{E}\left(
\frac{\|S_{n}\|}{n^{1/p}}\right)^{q} < \infty
\end{equation}
if and only if
\begin{equation}
X \in SLLN(p, q)
\end{equation}
and
\begin{equation}
\left \{
\begin{array}{ll}
\mbox{$\displaystyle \int_{0}^{\infty}
\mathbb{P}^{q/p}\left(\|X\|^{q}
> t \right) dt < \infty$} & \mbox{if $0 < q < p$,}\\
&\\
\mbox{$\displaystyle \mathbb{E}\|X\|^{p} \ln (1 + \|X\|) < \infty$}
& \mbox{if $q = p$,}\\
&\\
\mbox{$\displaystyle \mathbb{E}\|X\|^{q} < \infty$}
& \mbox{if $q > p$.}
\end{array}
\right.
\end{equation}
Furthermore, each of (1.1) and (1.2) implies that
\begin{equation}
\lim_{n \rightarrow \infty} \frac{S_{n}}{n^{1/p}} = 0~~{a.s.}
\end{equation}
For $0 < q < p$, (1.1) and (1.2) are equivalent so that each of them
implies that (1.3) and (1.4) hold.
\end{theorem}

\vskip 0.3cm

\begin{remark}
Let $q = 1$. Then one can easily see that Theorems 2.1 and 2.2 of Li, Qi, and Rosalsky [10]
follow from Theorem 1.4.
\end{remark}

\vskip 0.3cm

\begin{remark}
It follows from the conclusion (1.4) of Theorem 1.4 that, if (1.2) holds for some $q = q_{1} > 0$
then (1.2) holds for all $q > q_{1}$.
\end{remark}

\vskip 0.3cm

The current work continues the investigations by Hechner and Heinkel [5]
and Li, Qi, and Rosalsky [10] and [11]. More specifically:

\begin{description}
\item
{\bf (i)} For $0 < p < 1$ and $p < q < \infty$ and without any conditions
being imposed on the Banach space $\mathbf{B}$ we obtain in Theorem 2.1
necessary and sufficient conditions for $X \in SLLN(p, q)$.

\item
{\bf (ii)} For $1 \leq q < \infty$ we obtain assuming the Banach space
$\mathbf{B}$ is of stable type $p$ where $1 < p < 2$ (Theorem 2.2) or $p = 1$
(Theorem 2.3) necessary and sufficient conditions for $X \in SLLN(p, q)$.
\end{description}

Theorems 1.4, 2.1, 2.2, and 2.3 are new results when $\mathbf{B} = \mathbb{R}$
(Theorem 2.4).

When $\mathbf{B} = \mathbb{R}$, necessary and sufficient conditions
for $X \in SLLN(p, q)$ for the case where $0 < q < 1 \leq p < 2$ and for the case where
$0 < q \leq p < 1$ remain open problems.

The plan of the paper is as follows. The main results are stated in
Section 2 and they are proved in Section 3. In Section 4, three examples
will be provided for illustrating the necessary and sufficient conditions
obtained in this paper.

\section{Statement of the main results}

With the preliminaries accounted for, the main results may be stated.

\vskip 0.3cm

\begin{theorem}
Let $0 < p < 1$ and $p < q < \infty$.  Let $ \{X_{n};~n \geq 1 \}$ be a sequence
of independent copies of a $\mathbf{B}$-valued random variable $X$. Then we have
the following two statements:
\begin{description}
\item {\bf (a)}~~$X \in SLLN(p, q)$ if and only if
$\displaystyle \mathbb{E} \|X\|^{p} < \infty$,

\item {\bf (b)}~~$\displaystyle \sum_{n = 1}^{\infty} \frac{1}{n} \mathbb{E}\left(
\frac{\|S_{n}\|}{n^{1/p}}\right)^{q} < \infty$ if and only if
$\displaystyle \mathbb{E}\|X\|^{q} < \infty$.
\end{description}
\end{theorem}

\vskip 0.3cm

Let $X$ be a {\bf B}-valued random variable. For each $n \geq 1$, we
define the {\it quantile} $u_{n}$ of order $1 - \frac{1}{n}$
of $\|X\|$ as follows:
\[
u_{n} = \inf\left\{t:~ \mathbb{P}(\|X\| \leq t) > 1 - \frac{1}{n}
\right\} = \inf \left\{t:~ \mathbb{P}(\|X\| > t) < \frac{1}{n}
\right\}.
\]
If $\mathbb{E}\|X\| < \infty$, then it is easy to show that
\[
\lim_{n \rightarrow \infty} \frac{u_{n}}{n} = 0.
\]

\vskip 0.3cm

\begin{theorem}
Let $1 < p < 2$ and $1 \leq q < \infty$. Let $\mathbf{B}$ be a Banach space of stable type $p$. Let
$\{X_{n};~n \geq 1 \}$ be a sequence of independent copies of a $\mathbf{B}$-valued random variable $X$.
Then
\begin{equation}
X \in SLLN(p, q)
\end{equation}
if and only if
\begin{equation}
\left \{
\begin{array}{l}
\mbox{$\displaystyle \mathbb{E}X = 0$~~and}\\
~\\
\left\{
\begin{array}{ll}
\mbox{$\displaystyle \int_{0}^{\infty}
\mathbb{P}^{q/p}\left(\|X\|^{q}
> t \right) dt < \infty$} & \mbox{if $1 \leq q < p$,}\\
&\\
\mbox{$\displaystyle \mathbb{E}\|X\|^{p} < \infty
~~\mbox{and}~~\sum_{n=1}^{\infty} \frac{\int_{\min\{u_{n}^{p},
n\}}^{n} \mathbb{P}\left(\|X\|^{p} > t \right) dt}{n} < \infty$}
& \mbox{if~~$\displaystyle q = p $,}\\
&\\
\mbox{$\displaystyle \mathbb{E}\|X\|^{p} < \infty$}
& \mbox{if $q > p$.}\\
\end{array}
\right.
\end{array}
\right.
\end{equation}
\end{theorem}

\vskip 0.3cm

\begin{remark}
When $q = 1$ and $\mathbf{B}$ is of stable $p$ where $1 < p < 2$, Corollary 2.1 of Li, Qi,
and Rosalsky [10] follows immediately from Theorems 1.4 and 2.2; that is, (1.1),
(1.2), and (2.2) are equivalent.
\end{remark}

\vskip 0.3cm

Note by Lemma 5.6 of Li, Qi, and Rosalsky [10] that
\[
\sum_{n=1}^{\infty} \frac{\int_{\min\{u_{n}^{p}, n\}}^{n}
\mathbb{P}\left(\|X\|^{p} > t \right) dt}{n} < \infty
~~\mbox{whenever}~~\mathbb{E}\|X\|^{p} \ln^{\delta}(1 + \|X\|) < \infty
~~\mbox{for some}~~\delta > 0.
\]
Thus, for the interesting case $q = p$, Theorem 2.2 yields the following result.

\vskip 0.3cm

\begin{corollary}
Let $1 < p < 2$ and let $\{X_{n};~n \geq 1 \}$ be a sequence of independent
copies of a $\mathbf{B}$-valued random variable $X$.
If $\mathbf{B}$ is of stable type $p$, then
\[
X \in SLLN(p, p) \mbox{~whenever}~~\mathbb{E}X = 0 ~~\mbox{and}
~~\mathbb{E}\|X\|^{p} \ln^{\delta}(1 + \|X\|) < \infty
~~\mbox{for some}~~\delta > 0.
\]
\end{corollary}

\vskip 0.3cm

For the case where $1 < p < 2$ and $1 \leq q < \infty$, combining
Theorems 1.4 and 2.2, we immediately obtain necessary and sufficient
conditions for (1.1) to hold assuming that $\mathbf{B}$ is of stable type $p$.

\vskip 0.3cm
\begin{corollary}
Let $1 < p < 2$ and $1 \leq q < \infty$. Let $X$ be a $\mathbf{B}$-valued random
variable. If $\mathbf{B}$ is of stable type $p$, then (1.1) holds if and only if
\[
\left \{
\begin{array}{l}
\mbox{$\displaystyle \mathbb{E}X = 0$ ~ and}\\
~\\
\left \{
\begin{array}{ll}
\mbox{$\displaystyle \int_{0}^{\infty}
\mathbb{P}^{q/p}\left(\|X\|^{q}
> t \right) dt < \infty$} & \mbox{if $1 \leq q < p$,}\\
&\\
\mbox{$\displaystyle \mathbb{E}\|X\|^{p} \ln (1 + \|X\|) < \infty$}
& \mbox{if $q = p$,}\\
&\\
\mbox{$\displaystyle \mathbb{E}\|X\|^{q} < \infty$}
& \mbox{if $q > p$.}
\end{array}
\right.
\end{array}
\right.
\]
\end{corollary}

\vskip 0.3cm

\begin{remark}
For the case where $q = 1$, Corollary 2.2 above is Theorem 1.3 (i.e.,
Theorem 5 of Hechner and Heinkel [5]). Actually Corollary 2.2 for the
case where $q = 1$ is somewhat stronger than Theorem 5 (necessity half)
of Hechner and Heinkel [5] because $\mathbb{E}X = 0$ is an assumption
in Theorem 5 of Hechner and Heinkel [5].
\end{remark}

\vskip 0.3cm

We now present necessary and sufficient conditions for (1.2) for the
case where $p = 1$ and $1 \leq q < \infty$.

\vskip 0.3cm

\begin{theorem}
Let $1 \leq q < \infty$ and let $\mathbf{B}$ be a Banach space of stable
type $1$. Let $\{X_{n};~n \geq 1 \}$ be a sequence of independent copies
of a $\mathbf{B}$-valued random variable $X$. Then
\begin{equation}
X \in SLLN(1, q)
\end{equation}
if and only if
\begin{equation}
\left\{
\begin{array}{l}
\mbox{$\displaystyle \mathbb{E}\|X\| < \infty, ~\mathbb{E}X = 0$,~~and}\\
~\\
\left\{
\begin{array}{ll}
\mbox{$\displaystyle
\sum_{n=1}^{\infty} \frac{\left\|\mathbb{E}XI\{\|X\| \leq n\}
\right\|}{n} < \infty$~~and~~$\displaystyle
\sum_{n=1}^{\infty} \frac{\int_{\min\{u_{n}, n\}}^{n}
\mathbb{P}(\|X\| > t) dt}{n} < \infty$} & \mbox{if $q = 1$,}\\
&\\
\mbox{$\displaystyle
\sum_{n=1}^{\infty} \frac{\left\|\mathbb{E}XI\{\|X\| \leq n\}
\right\|^{q}}{n} < \infty$} & \mbox{if $q > 1$.}
\end{array}
\right.
\end{array}
\right.
\end{equation}
\end{theorem}

\vskip 0.3cm

\begin{remark}
For the case where $q = 1$, Theorem 2.3 is Theorem 2.3 of
Li, Qi, and Rosalsky [10].
\end{remark}

\vskip 0.3cm

By Lemmas 5.5 and 5.6 of Li, Qi, and Rosalsky [10], (2.4) holds whenever $\mathbb{E}X = 0$ and
$\mathbb{E}\|X\|\ln(1 + \|X\|) < \infty$. Combining Theorems 1.4 and 2.3, we immediately have
the following result.

\vskip 0.3cm

\begin{corollary}
Let $1 \leq q < \infty$ and let $\mathbf{B}$ be a Banach space of stable
type $1$. Let $\{X_{n};~n \geq 1 \}$ be a sequence of independent copies
of a $\mathbf{B}$-valued random variable $X$. Then
\[
\sum_{n=1}^{\infty} \frac{1}{n} \mathbb{E}\left(\frac{\|S_{n}\|}{n} \right)^{q}
< \infty
\]
if and only if
\[
\left\{
\begin{array}{l}
\mbox{$\displaystyle \mathbb{E}X = 0$~~and}\\
~\\
\left\{
\begin{array}{ll}
\mbox{$\displaystyle
\mathbb{E}\|X\| \ln (1+\|X\|) < \infty$} & \mbox{if $q = 1$,}\\
&\\
\mbox{$\displaystyle
\mathbb{E} \|X\|^{q} < \infty$}
& \mbox{if $q > 1$.}
\end{array}
\right.
\end{array}
\right.
\]
\end{corollary}

\vskip 0.3cm

As a summary of our Theorems 1.4 and 2.1-2.3 and Corollaries 2.2 and 2.3, we now
present the following theorem for a real-valued random variable $X$. For
$q = 1$, the equivalence of (i) and (ii) has recently been obtained by
Li, Qi, and Rosalsky [10], and for $1 = q < p < 2$, the equivalence
of (iii) and (iv) is due to Hechner and Heinkel [5] (see Theorem 1.3 above)
assuming that $\mathbb{E}X = 0$ for the implication ((iii) $\Rightarrow$ (iv)).

\vskip 0.3cm

\begin{theorem}
Let $0 < p < 2$ and $1 \leq q < \infty$. Let $\{X_{n}; ~n \geq 1 \}$ be a
sequence of independent copies of a real-valued random variable $X$. The
following two statements are equivalent:
\begin{description}
\item {\bf (i)} \quad
$\displaystyle X \in SLLN(p, q)$,

\item {\bf (ii)} \quad
$\displaystyle
\left \{
\begin{array}{ll}
\mbox{$\displaystyle \mathbb{E}X = 0 ~~\mbox{and}~~\int_{0}^{\infty}
\mathbb{P}^{q/p}\left(|X|^{q}
> t \right) dt < \infty$} & \mbox{if $1 \leq q < p < 2$,}\\
&\\
\mbox{$\displaystyle \mathbb{E}X = 0, ~~\mathbb{E} |X|^{p} < \infty, ~~\mbox{and}
~~\sum_{n=1}^{\infty} \frac{\int_{\min\{u_{n}^{p},
n\}}^{n} \mathbb{P}\left(|X|^{p} > t \right) dt}{n} < \infty$} & \mbox{if
~~$\displaystyle 1 < q = p < 2$,}\\
&\\
\mbox{$\displaystyle \mathbb{E}X = 0 ~~\mbox{and}~~\mathbb{E}|X|^{p} < \infty$}
& \mbox{if $1 < p < 2$ \mbox{and} $q > p$,}\\
&\\
\mbox{$\displaystyle \mathbb{E}X = 0, ~~\sum_{n=1}^{\infty}
\frac{|\mathbb{E}XI\left\{|X| \leq n \right\}|}{n} < \infty,~\mbox{and}$} &\\
&\\
\mbox{$\displaystyle \sum_{n=1}^{\infty} \frac{\int_{\min\{u_{n}, n\}}^{n}
\mathbb{P}\left(|X| > t \right) dt}{n} < \infty$} & \mbox{if
~~$\displaystyle q = p = 1$,}\\
&\\
\mbox{$\displaystyle \mathbb{E}X = 0~~\mbox{and}~~\sum_{n=1}^{\infty}
\frac{\left|\mathbb{E}XI\left\{|X| \leq n \right\} \right|^{q}}{n} < \infty$}
& \mbox{if $p = 1 < q$,}\\
&\\
\mbox{$\displaystyle \mathbb{E}|X|^{p} < \infty$}
& \mbox{if $0 < p < 1 \leq q$.}
\end{array}
\right.
$
\end{description}
The following two statements are equivalent:
\begin{description}
\item {\bf (iii)} \quad
$\displaystyle \sum_{n = 1}^{\infty} \frac{1}{n}\mathbb{E} \left(\frac{\left|S_{n}\right|}{n^{1/p}}
\right)^{q} < \infty,$

\item {\bf (iv)} \quad
$\displaystyle
\left \{
\begin{array}{ll}
\mbox{$\displaystyle \mathbb{E}X = 0 ~~\mbox{and}~~\int_{0}^{\infty}
\mathbb{P}^{q/p}\left(|X|^{q}
> t \right) dt < \infty$} & \mbox{if $1 \leq q < p < 2$,}\\
&\\
\mbox{$\displaystyle \mathbb{E}X = 0 ~~\mbox{and}~~\mathbb{E} |X|^{p}\ln (1 + |X|) < \infty$}
& \mbox{if~~$\displaystyle 1 \leq q = p < 2$,}\\
&\\
\mbox{$\displaystyle \mathbb{E}X = 0 ~~\mbox{and}~~\mathbb{E}|X|^{p} < \infty$}
& \mbox{if $1 \leq p < 2$ \mbox{and} $q > p$,}\\
&\\
\mbox{$\displaystyle \mathbb{E}|X|^{q} < \infty$}
& \mbox{if $0 < p < 1 \leq q$.}
\end{array}
\right.
$
\end{description}
\end{theorem}

\section{Proofs of Theorems 2.1 - 2.3}

In this section we denote by $C_{k}$ positive constants the precise values of which
do not matter.

First we introduce some notation. Let $(a_{k})_{1 \leq k \leq n}$ be a finite
sequence of real numbers and $(a^{*}_{k})_{1 \leq k \leq n}$ the nonincreasing
rearrangement of the sequence $(|a_{k}|)_{1 \leq k \leq n}$. For a given $r \geq 1$,
\[
\left\|(a_{k})_{1 \leq k \leq n} \right\|_{r, \infty} = \sup_{1 \leq k \leq n} k^{1/r} a^{*}_{k}
\]
is called the {\it weak}-$\ell_{r}$ {\it norm} of the sequence $(a_{k})_{1 \leq k \leq n}$.
Let $V_{k}, 1 \leq k \leq n$ be independent real-valued random variables. Then the remarkable
Marcus-Pisier [14] inequality asserts that for all $r \geq 1$,
\begin{equation}
\mathbb{P} \left(\left\|(V_{k})_{1 \leq k \leq n} \right\|_{r, \infty} > u \right)
\leq \frac{2e}{u^{r}} \sup_{t > 0} \left(t^{r} \sum_{k=1}^{n}
\mathbb{P} \left(|V_{k}| > t \right) \right) ~~\forall ~u > 0.
\end{equation}
The original Marcus-Pisier [14] inequality involved the constant $262$ instead of $2e$.
The improved constant is due to J. Zinn (see Pisier [18, Lemma 4.11]).

Let $X$ be a {\bf B}-valued random variable. For each $n \geq 1$, let the quantile $u_{n}$ of order
$1 - \frac{1}{n}$ of $\|X\|$ be defined as in Section 2. We then see that for every $q > 0$,
\[
\inf\left\{t:~ \mathbb{P}\left(\|X\|^{q} \leq t \right) > 1 - \frac{1}{n}
\right\} = \inf \left\{t:~ \mathbb{P}\left(\|X\|^{q} > t \right) < \frac{1}{n}
\right\} = u_{n}^{q};
\]
i.e., $u_{n}^{q}$ is the quantile of order $1 - \frac{1}{n}$
of $\|X\|^{q}$. Let $\{X_{n}; ~n \geq 1 \}$ be a sequence of independent
copies of $\mathbf{B}$-valued variable $X$. Write, for $n \geq 1$,
\[
S^{(1)}_{n} = \sum_{k=1}^{n} X_{k} I\{\|X_{k}\|^{p} \leq k\},~~S^{(2)}_{n}
= S_{n} - S^{(1)}_{n} = \sum_{k=1}^{n} X_{k} I\{\|X_{k}\|^{p} > k\},
\]
\[
U_{n} = \sum_{k=1}^{n} X_{k} I\{\|X_{k}\|^{p} \leq n \}, ~~
U_{n}^{(1)} = \sum_{k=1}^{n} X_{k}I\{\|X_{k}\| \leq u_{n} \},
~~\mbox{and}~~U_{n}^{(2)} = U_{n} - U_{n}^{(1)},~~n \geq 1.
\]
Motivated by Lemma 1 of Hechner and Heinkel [5] and its proof, we establish
the following result.

\vskip 0.3cm

\begin{lemma}
Let $1 < p < 2$ and $1 \leq q < p$. Let $\mathbf{B}$ be a Banach space of stable type $p$. Then
there exists a universal constant $c(p, q) > 0$ such that, for every finite sequence
$V_{k}, 1 \leq k \leq n$ of independent $\mathbf{B}$-valued random variables
with $\max_{1\leq k \leq n} \mathbb{E}\|V_{k}\|^{q} < \infty$,
\begin{equation}
\mathbb{E}\left\|\sum_{k=1}^{n}\left(V_{k} - \mathbb{E}V_{k}\right) \right\|^{q} \leq c(p, q)
\left( \sup_{t > 0} t^{p/q} \sum_{k=1}^{n}
\mathbb{P} \left(\|V_{k}\|^{q} > t \right) \right)^{q/p}.
\end{equation}
\end{lemma}

\vskip 0.3cm

\begin{remark}
Clearly, if $q = 1$, then Lemma 3.1 is Lemma 1 of Hechner and Heinkel [5].
\end{remark}

{\it Proof of Lemma 3.1}~~Let $\{V^{\prime}_{k};~1 \leq k \leq n\}$
be an independent copy of $\{V_{k};~1 \leq k \leq n\}$
and let $\{R_{k}; ~1 \leq k \leq n\}$ be a Rademacher sequence independent of
$\{V_{k}, V^{\prime}_{k};~1 \leq k \leq n\}$. Since
$q \geq 1$, $g(x) = x^{p}, ~x \in [0, \infty)$ is a convex nonnegative function.
Applying (2.5) of Ledoux and Talagrand [9, p. 46], we have that
\begin{equation}
\mathbb{E}\left\|\sum_{k=1}^{n}\left(V_{k} - \mathbb{E}V_{k}\right) \right\|^{q}
\leq \mathbb{E} \left\|\sum_{k=1}^{n}\left(V_{k} - V^{\prime}_{k}\right)\right\|^{q}
= \mathbb{E} \left\| \sum_{k=1}^{n}R_{k}\left(V_{k} - V^{\prime}_{k}\right)\right\|^{q}
\leq 2^{q-1} \mathbb{E} \left\| \sum_{k=1}^{n}R_{k}V_{k}\right\|^{q}.
\end{equation}
Since $\mathbf{B}$ is of stable type $p$
with $1 \leq p < 2$, the Maurey-Pisier [16] theorem
asserts that it is also of stable type $r$ for some $r > p$.
Let $\left(A^{*}_{k} \right)_{1 \leq k \leq n}$
be the nonincreasing rearrangement of $\left(\|V_{k}\| \right)_{1 \leq k \leq n}$.
Note that $r/q > 1$, $p/q > 1$ (since $1 \leq q < p < r$), and $\mathbf{B}$ is also of
Rademacher type $r$. We thus have that
\begin{equation}
\begin{array}{lll}
\mbox{$\displaystyle \mathbb{E} \left\| \sum_{k=1}^{n}R_{k}V_{k}\right\|^{q}$}
&=& \mbox{$\displaystyle
\mathbb{E} \left(\mathbb{E}\left( \left.
\left\| \sum_{k=1}^{n}R_{k}V_{k}\right\|^{q} \right|V_{1}, ..., V_{n} \right) \right)$}\\
&&\\
&\leq& \mbox{$\displaystyle
\mathbb{E} \left(\mathbb{E}\left( \left. \left\| \sum_{k=1}^{n}R_{k}V_{k}\right\|^{r}
\right|V_{1}, ..., V_{n} \right) \right)^{1/(r/q)}$}\\
&&\\
&\leq& \mbox{$\displaystyle C_{1}
\mathbb{E} \left(\sum_{k=1}^{n} \|V_{k}\|^{r} \right)^{q/r}$}\\
&&\\
&=& \mbox{$\displaystyle C_{1}
\mathbb{E} \left( \sum_{k=1}^{n}\left(k^{r/p}
\left(A^{*}_{k} \right)^{r} \right) k^{-r/p} \right)^{q/r}$}\\
&&\\
&\leq& \mbox{$\displaystyle C_{1}
\mathbb{E} \left(\left(\sup_{1 \leq k \leq n} k^{q/p} (A^{*}_{k})^{q} \right)
\left(\sum_{k=1}^{n} k^{-r/p} \right)^{q/r} \right)$}\\
&&\\
&=& \mbox{$\displaystyle C_{2} \mathbb{E}
\left\|\left(\|V_{k}\|^{q} \right)_{1 \leq k \leq n} \right\|_{p/q, \infty}.$}
\end{array}
\end{equation}
Write $\Delta = \sup_{t > 0} t^{p/q}
\sum_{k=1}^{n} \mathbb{P} \left(\|V_{k}\|^{q} > t \right)$.
Using the Marcus-Pisier [14] inequality (3.1), we have that
\begin{equation}
\begin{array}{lll}
\mbox{$\displaystyle \mathbb{E}
\left\|\left(\|V_{k}\|^{q} \right)_{1 \leq k \leq n} \right\|_{p/q, \infty}$}
&=&
\mbox{$\displaystyle \left(\int_{0}^{\Delta^{q/p}} + \int_{\Delta^{q/p}}^{\infty} \right)
\mathbb{P}\left(
\left\|\left(\|V_{k}\|^{q} \right)_{1 \leq k \leq n} \right\|_{p/q, \infty} > t \right)dt$}\\
&&\\
&\leq& \mbox{$\displaystyle \Delta^{q/p}
+ \int_{\Delta^{q/p}}^{\infty} \frac{2e \Delta}{t^{p/q}}dt$}\\
&&\\
&=& \mbox{$\displaystyle \left(1 + \frac{2qe}{p-q} \right)\Delta^{q/p}.$}
\end{array}
\end{equation}
Now (3.2) follows from (3.3), (3.4), and (3.5).~ $\Box$

\vskip 0.3cm

The following nice result is Proposition 3 of Hechner and Heinkel [5].

\vskip 0.3cm

\begin{lemma}
{\rm (Hechner and Heinkel [5])}. Let $p > 1$ and let $ \{X_{n};~n \geq 1 \}$
be a sequence of independent copies of a $\mathbf{B}$-valued random variable $X$.
Write
\[
u_{n} = \inf \left\{t:~ \mathbb{P}(\|X\| > t) < \frac{1}{n}
\right\}, ~n \geq 1.
\]
Then the following three statements are equivalent:
\begin{description}
\item
\quad {\bf (i)} \quad $\displaystyle \int_{0}^{\infty}
\mathbb{P}^{1/p}(\|X\| > t) dt < \infty;$

\item
\quad {\bf (ii)} \quad $\displaystyle \sum_{n=1}^{\infty}
\frac{u_{n}}{n^{1 + 1/p}} < \infty;$

\item
\quad {\bf (iii)} \quad
$\displaystyle \sum_{n=1}^{\infty} \frac{1}{n^{1 + 1/p}}
\mathbb{E} \left(\max_{1 \leq k \leq n} \|X_{k}\| \right) < \infty.$
\end{description}
\end{lemma}

\vskip 0.3cm

The next lemma and its proof are similar to Lemma 3 of Hechner and Heinkel [5]
and its proof, respectively.

\vskip 0.3cm

\begin{lemma}
Let $1 \leq q < p < 2$. Let $X$ be a {\bf B}-valued random variable with
\begin{equation}
\int_{0}^{\infty} \mathbb{P}^{q/p}\left(\|X\|^{q} > t \right) dt < \infty.
\end{equation}
If $\mathbf{B}$ is a Banach space of Rademacher type $q$, then
\begin{equation}
\sum_{n=1}^{\infty} \frac{\mathbb{E}
\left\|\left(S_{n} - U_{n}^{(1)} \right)
- \mathbb{E}\left(S_{n} - U_{n}^{(1)} \right)\right\|^{q}}
{n^{1 + q/p}} < \infty.
\end{equation}
\end{lemma}

{\it Proof}~~Let $f_{q}(t) = \mathbb{P}\left(\|X\|^{q} > t \right), ~t \geq 0$.
Since $\mathbf{B}$ is a Banach space of Rademacher type $q$ and
\[
\left(S_{n} - U_{n}^{(1)} \right) - \mathbb{E}\left(S_{n} - U_{n}^{(1)} \right)
= \sum_{k=1}^{n}\left(X_{k}I\{\|X_{k}\| > u_{n}\}
- \mathbb{E}XI\{\|X\| > u_{n}\}\right), ~~n \geq 1,
\]
we have that
\begin{equation}
\begin{array}{lll}
\mbox{$\displaystyle
\mathbb{E} \left\|\left(S_{n} - U_{n}^{(1)} \right)
- \mathbb{E}\left(S_{n} - U_{n}^{(1)} \right)\right\|^{q} $}
&\leq & \mbox{$\displaystyle C_{3} n \mathbb{E} \left\|XI\{\|X\| > u_{n}\}
- \mathbb{E}XI\{\|X\| > u_{n}\} \right\|^{q}$}\\
&&\\
&\leq& \mbox{$\displaystyle
C_{4} n \mathbb{E} \left(\|X\|^{q} I\left\{\|X\|^{q} > u_{n}^{q} \right\}\right)$}\\
&&\\
&=& \mbox{$\displaystyle
C_{4} \left(nu_{n}^{p} \mathbb{P}\left(\|X\|^{q} > u_{n}^{q} \right)
+ n \int_{u_{n}^{q}}^{\infty} f_{q}(t)dt \right)$}\\
&&\\
&\leq& \mbox{$\displaystyle C_{4}
\left(u_{n}^{q} + n \int_{u_{n}^{q}}^{\infty} f_{q}(t)dt \right).$}\\
\end{array}
\end{equation}
Set $p_{1} = p/q$, $Y = \|X\|^{q}$, and $u_{n,q} = u_{n}^{q}, ~n \geq 1$. Noting that
$p_{1} > 1$ (since $1 \leq q < p < 2$), by Lemma 3.2 (i.e., Proposition 3 of Hechner
and Heinkel [5]), it follows from (3.6) that
\begin{equation}
\sum_{n=1}^{\infty} \frac{u_{n}^{q}}{n^{1 + q/p}}
= \sum_{n=1}^{\infty} \frac{u_{n,q}}{n^{1 + 1/p_{1}}} < \infty.
\end{equation}
Also (3.6) implies that
\begin{equation}
\begin{array}{lll}
\mbox{$\displaystyle
\sum_{n=1}^{\infty} \frac{1}{n^{q/p}} \int_{u_{n}^{q}}^{\infty} f_{q}(t)dt$}
&=&
\mbox{$\displaystyle
\sum_{n=1}^{\infty} \frac{1}{n^{1/p_{1}}} \int_{u_{n,q}}^{\infty} f_{q}(t)dt$}\\
&&\\
&=&
\mbox{$\displaystyle
\sum_{n=1}^{\infty} n^{-1/p_{1}} \sum_{j=n}^{\infty}
\int_{u_{j,q}}^{u_{j+1, q}} f_{q}(t)dt$}\\
&&\\
&=&
\mbox{$\displaystyle
\sum_{j=1}^{\infty} \left(\int_{u_{j,q}}^{u_{j+1, q}} f_{q}(t)dt \right)
\sum_{n=1}^{j} n^{-1/p_{1}}$}\\
&&\\
&\leq&
\mbox{$\displaystyle
C_{5}\sum_{j=1}^{\infty}
\left(\int_{u_{j,q}}^{u_{j+1, q}} f_{q}^{1/p_{1}}(t)dt \right)
\frac{j^{1- 1/p_{1}}}{j^{1- 1/p_{1}}}$}\\
&&\\
&\leq&
\mbox{$\displaystyle
C_{5} \int_{0}^{\infty} f_{q}^{1/p_{1}}(t)dt$}\\
&&\\
&=&
\mbox{$\displaystyle C_{5} \int_{0}^{\infty}
\mathbb{P}^{q/p}\left(\|X\|^{q} > t \right)dt$}\\
&&\\
&<& \mbox{$\displaystyle \infty.$}
\end{array}
\end{equation}
The conclusion (3.7) follows from (3.8), (3.9), and (3.10). ~$\Box$

\vskip 0.3cm

The proof of the next lemma is similar to that of Lemma 4 of Hechner
and Heinkel [5] and Lemma 5.3 of Li, Qi, and Rosalsky [10] and it 
contains a nice application of Lemma 3.1 above.

\vskip 0.3cm

\begin{lemma}
Let $1 \leq q \leq p < 2$. Let $X$ be a {\bf B}-valued random variable with (3.6).
If $\mathbf{B}$ is a Banach space of stable type $p$, then
\begin{equation}
\sum_{n=1}^{\infty} \frac{\mathbb{E} \left\|U_{n}^{(1)} - \mathbb{E}U_{n}^{(1)} \right\|^{q}}
{n^{1 + q/p}} < \infty.
\end{equation}
\end{lemma}

\vskip 0.3cm

\begin{remark}
Note that
\[
\mathbb{E}\|X\|^{q} = \int_{0}^{\infty} \mathbb{P} \left(\|X\|^{q} > t \right)dt.
\]
Thus for $q = p$, (3.6) holds if and only if $\mathbb{E}\|X\|^{p} < \infty$. By Lemma 3.4, if
$\mathbf{B}$ is a Banach space of stable type $p \in [1, 2)$, then
\begin{equation}
\sum_{n=1}^{\infty} \frac{\mathbb{E} \left\|U_{n}^{(1)} - \mathbb{E}U_{n}^{(1)} \right\|^{p}}
{n^{2}} < \infty
\end{equation}
whenever $\mathbb{E}\|X\|^{p} < \infty$.
\end{remark}

\vskip 0.3cm

{\it Proof of Lemma 3.4}~~Since $\mathbf{B}$ is of stable type $p$, the Maurey-Pisier [16]
theorem asserts that it is also of stable type $r$ for some $r > p$. Applying Lemma 3.1,
there exists a universal constant $0 < c(r, q) < \infty$ such that
\[
\begin{array}{lll}
\mbox{$\displaystyle
\mathbb{E} \left\|U_{n}^{(1)} - \mathbb{E}U_{n}^{(1)} \right\|^{q}$}
&\leq&
\mbox{$\displaystyle
c(r, q) \left(\sup_{t > 0} t^{r/q}
\sum_{k=1}^{n}
\mathbb{P} \left(\|X_{k}\|^{q}I\{\|X_{k}\| \leq u_{n} \} > t \right)\right)^{q/r}$}\\
&&\\
&\leq&
\mbox{$\displaystyle
c(r, q) \left(n \sup_{0 \leq t \leq u_{n}^{q}}
t^{r/q} \mathbb{P}\left(\|X\|^{q} > t \right) \right)^{q/r}, ~~n \geq 1.$}
\end{array}
\]
It is easy to see that for all $x > 0$,
\[
\begin{array}{lll}
\mbox{$\displaystyle
\left(\int_{0}^{x} \mathbb{P}^{q/r} \left(\|X\|^{q} > t \right)dt \right)^{r/q}$}
&\geq&
\mbox{$\displaystyle
\left(\int_{0}^{x} \mathbb{P}^{q/r} \left(\|X\|^{q} > x \right)dt \right)^{r/q}$}\\
&&\\
&=&
\mbox{$\displaystyle
x^{r/q} \mathbb{P} \left(\|X\|^{q} > x \right).$}
\end{array}
\]
We thus have that
\[
\begin{array}{lll}
\mbox{$\displaystyle
\mathbb{E} \left\|U_{n}^{(1)} - \mathbb{E}U_{n}^{(1)} \right\|^{q}$}
&\leq&
\mbox{$\displaystyle
c(r, q)  \left(n \sup_{0 \leq t \leq u_{n}^{q}}
t^{r/q} \mathbb{P}\left(\|X\|^{q} > t \right) \right)^{q/r}$}\\
&&\\
&\leq&
\mbox{$\displaystyle
c(r, q) n^{q/r} \int_{0}^{u_{n}^{q}} \mathbb{P}^{q/r}\left(\|X\|^{q} > t \right)dt, ~~n \geq 1.$}
\end{array}
\]
Let $u_{0} = 0$ and note that $\mathbb{P}\left(\|X\|^{q} > t \right) \geq 1/k$ for
$t \in [u_{k-1}^{q}, u_{k}^{q}), ~k \geq 1$. It follows that
\[
\begin{array}{lll}
\mbox{$\displaystyle
\sum_{n=1}^{\infty}
\frac{\mathbb{E}\left\|U_{n}^{(1)} - \mathbb{E}U_{n}^{(1)}
\right\|^{q}}{n^{1 + q/p}} $}
&\leq& \mbox{$\displaystyle c(r, q) \sum_{n=1}^{\infty} \frac{1}{n^{1 + q/p -q/r}}
\int_{0}^{u_{n}^{q}}
\mathbb{P}^{q/r}\left(\|X\|^{q} > t \right) dt$}\\
&&\\
&=& \mbox{$\displaystyle c(r,q)\sum_{n=1}^{\infty} \frac{1}{n^{1 + q/p -q/r}}
\sum_{k=1}^{n}\int_{u_{k-1}^{q}}^{u_{k}^{q}} \mathbb{P}^{q/r}\left(\|X\|^{q}
> t \right) dt$}\\
&&\\
&=& \mbox{$\displaystyle c(r,q)\sum_{k=1}^{\infty}
\left(\sum_{n=k}^{\infty} \frac{1}{n^{1 + q/p -q/r}} \right)
\int_{u_{k-1}^{q}}^{u_{k}^{q}} \mathbb{P}^{q/r}\left(\|X\|^{q} > t \right) dt$}\\
&&\\
&\leq& \mbox{$\displaystyle C_{6} \sum_{k=1}^{\infty}\frac{1}{k^{q/p -q/r}}
\int_{u_{k-1}^{q}}^{u_{k}^{q}} \mathbb{P}^{q/r}\left(\|X\|^{q} > t \right) dt$}\\
&&\\
&\leq& \mbox{$\displaystyle C_{6} \sum_{k=1}^{\infty}
\int_{u_{k-1}^{q}}^{u_{k}^{q}} \mathbb{P}^{q/p}\left(\|X\|^{q} > t \right) dt$}\\
&&\\
&=& \mbox{$\displaystyle
C_{6} \int_{0}^{\infty}\mathbb{P}^{q/p}\left(\|X\|^{q} > t \right) dt < \infty$}
\end{array}
\]
proving (3.11) and completing the proof of Lemma 3.4.~$\Box$

\vskip 0.3cm

\begin{lemma}
Let $1 \leq p < 2$ and let $X$ be a $\mathbf{B}$-valued random variable with
$\mathbb{E}\|X\|^{p} < \infty$. Then
\begin{equation}
\sum_{n=1}^{\infty} \frac{1}{n^{2}}
\left(\sum_{k=1}^{n} \mathbb{E} \|X\| I\left\{k < \|X\|^{p} \leq n \right\} \right)^{p}
< \infty,
\end{equation}
\begin{equation}
\sum_{n=1}^{\infty}
\frac{u_{n}^{p}}{n^{2}} < \infty,
\end{equation}
and for every $\delta > 0$,
\begin{equation}
\sum_{n=1}^{\infty}
\frac{\mathbb{E}\|X\|^{p+\delta} I\left\{\|X\|^{p} \leq n\right\}}{n^{1 + \delta/p}}
< \infty.
\end{equation}
Furthermore, if $p > 1$ then
\begin{equation}
\sum_{n=1}^{\infty}
\frac{\left(\mathbb{E}\|X\|I\left\{\|X\|^{p} > n\right\} \right)^{p}}{n^{2-p}}
< \infty.
\end{equation}
\end{lemma}

\vskip 0.3cm

\begin{remark}
For $p=1$, (3.13) and (3.14) together are Lemma 5.1 of Li, Qi, and Rosalsky [10].
\end{remark}

\vskip 0.3cm

{\it Proof of Lemma 3.5}~~Since $u_{n}^{p}$ is the quantile of order
$1 - \frac{1}{n}$ of $\|X\|^{p}$, (3.14) immediately follows
from the second half of Lemma 5.1 of Li, Qi, and Rosalsky [10].

The proof of (3.15) is easy and we leave it to the reader.

We now show that $\mathbb{E}\|X\|^{p} < \infty$ implies (3.13). For $n \geq 2$, let
\[
\Lambda_{n} = \sum_{k=2}^{n}k \mathbb{P}\left(k-1 < \|X\|^{p} \leq k \right),~~
\lambda_{n,j} = \frac{j\mathbb{P}\left(j-1 < \|X\|^{p} \leq j \right)}
{\Lambda_{n}}, ~~2 \leq j \leq n.
\]
Clearly
\[
\lambda_{n,j} \geq 0, ~2 \leq j \leq n, ~~\sum_{j=2}^{n} \lambda_{n,j} = 1,
\]
and
\[
\Lambda_{n} \leq \mathbb{E}\|X\|^{p} + 1 < \infty, ~~n \geq 2.
\]
Note that the function $\phi(t) = t^{p}$ is convex on $[0, \infty)$ and
\[
\begin{array}{lll}
\mbox{$\displaystyle
\|X\|\sum_{k=1}^{n}I\left\{k < \|X\|^{p} \leq n \right\}$}
&=&
\mbox{$\displaystyle
\sum_{k=1}^{n} \sum_{j=k+1}^{n} \|X\|I\left\{j-1 < \|X\|^{p} \leq j \right\}$}\\
&&\\
&\leq&
\mbox{$\displaystyle
\sum_{k=1}^{n} \sum_{j=k+1}^{n} j^{1/p}I\left\{j-1 < \|X\|^{p} \leq j \right\}$}\\
&&\\
&=&
\mbox{$\displaystyle
\sum_{j=2}^{n} \sum_{k=1}^{j-1} j^{1/p} I\left\{j-1 < \|X\|^{p} \leq j \right\}$}\\
&&\\
&\leq&
\mbox{$\displaystyle
\sum_{j=2}^{n} j^{1 + 1/p} I\left\{j-1 < \|X\|^{p} \leq j \right\}, ~~n \geq 2.$}
\end{array}
\]
We thus have that
\[
\begin{array}{lll}
\mbox{$\displaystyle
\left(\sum_{k=1}^{n} \mathbb{E} \|X\| I\left\{k < \|X\|^{p} \leq n \right\} \right)^{p}$}
&=&
\mbox{$\displaystyle
\left(\mathbb{E}\left(\|X\| \sum_{k=1}^{n} I\left\{k < \|X\|^{p} \leq n \right\} \right) \right)^{p}$}\\
&&\\
&\leq&
\mbox{$\displaystyle
\left(\mathbb{E}
\left(\sum_{j=2}^{n} j^{1 + 1/p} I\left\{j-1 < \|X\|^{p} \leq j \right\}\right) \right)^{p}$}\\
&&\\
&=&
\mbox{$\displaystyle
\left( \sum_{j=2}^{n} j^{1 + 1/p} \mathbb{P} \left(j-1 < \|X\|^{p} \leq j \right) \right)^{p}$}\\
&&\\
&=&
\mbox{$\displaystyle
\Lambda_{n}^{p} \left( \sum_{j=2}^{n} j^{1/p} \lambda_{n,j} \right)^{p}$}\\
&&\\
&\leq&
\mbox{$\displaystyle
\Lambda_{n}^{p} \sum_{j=2}^{n} \lambda_{n,j} \left(j^{1/p} \right)^{p} $}\\
&&\\
&=&
\mbox{$\displaystyle
\Lambda_{n}^{p-1} \sum_{j=2}^{n} j^{2}  \mathbb{P} \left(j-1 < \|X\|^{p} \leq j \right)$}\\
&&\\
&\leq&
\mbox{$\displaystyle
C_{7} \sum_{j=2}^{n} j^{2}  \mathbb{P} \left(j-1 < \|X\|^{p} \leq j \right)$}, ~~n \geq 2.
\end{array}
\]
It now is easy to see that
\[
\begin{array}{lll}
\mbox{$\displaystyle
\sum_{n=1}^{\infty} \frac{1}{n^{2}} \sum_{j=2}^{n} j^{2}
\mathbb{P} \left(j-1 < \|X\|^{p} \leq j \right)$}
&=&
\mbox{$\displaystyle
\sum_{j=2}^{\infty} \left(\sum_{n=j}^{\infty} \frac{1}{n^{2}} \right) j^{2}
\mathbb{P} \left(j-1 < \|X\|^{p} \leq j \right)$}\\
&&\\
&\leq&
\mbox{$\displaystyle
C_{8} \sum_{j=2}^{\infty} j \mathbb{P} \left(j-1 < \|X\|^{p} \leq j \right)$}\\
&&\\
&\leq&
\mbox{$\displaystyle
C_{8} \left(\mathbb{E}\|X\|^{p} + 1 \right) < \infty$}
\end{array}
\]
thereby proving (3.13).

We now prove (3.16). Note that for $n \geq 1$,
\[
\begin{array}{lll}
\mbox{$\displaystyle
\mathbb{E}\|X\|I\left\{\|X\|^{p} > n\right\}$}
&\leq&
\mbox{$\displaystyle \sum_{j=n+1}^{\infty}
j^{1/p} \mathbb{P} \left(j-1 < \|X\|^{p} \leq j \right)$}\\
&&\\
&=&
\mbox{$\displaystyle
\sum_{j=n+1}^{\infty} j^{1/p -1}
\left(j \mathbb{P} \left(j-1 < \|X\|^{p} \leq j \right)\right)$}
\end{array}
\]
and
\[
\sum_{j=n}^{\infty}
j \mathbb{P} \left(j-1 < \|X\|^{p} \leq j \right) \leq \mathbb{E}\|X\|^{p} + 1.
\]
Thus, by the same arguments used in proving (3.13), we have that
\[
\left(\mathbb{E}\|X\|I\left\{\|X\|^{p} > n\right\}\right)^{p}
\leq C_{9} \sum_{j=n}^{\infty} j^{1-p} \left(j \mathbb{P} \left(j-1 < \|X\|^{p} \leq j \right)\right),
~~n \geq 1.
\]
Since $p > 1$, we get that
\[
\begin{array}{lll}
\mbox{$\displaystyle
\sum_{n=1}^{\infty} \frac{\left(\mathbb{E}\|X\|I\left\{\|X\|^{p} > n\right\} \right)^{p}}{n^{2-p}}$}
&\leq&
\mbox{$\displaystyle
C_{9} \sum_{n=1}^{\infty} \frac{1}{n^{2-p}}
\sum_{j=n}^{\infty} j^{1-p} \left(j \mathbb{P} \left(j-1 < \|X\|^{p} \leq j \right)\right)$}\\
&&\\
&=&
\mbox{$\displaystyle
C_{9} \sum_{j=1}^{\infty} \left(\sum_{n=1}^{j} n^{p-2}\right)
j^{2-p} \mathbb{P} \left(j-1 < \|X\|^{p} \leq j \right)$}\\
&&\\
&\leq&
\mbox{$\displaystyle
C_{10} \sum_{j=1}^{\infty} j \mathbb{P} \left(j-1 < \|X\|^{p} \leq j \right)$}\\
&&\\
&\leq&
\mbox{$\displaystyle
C_{10} \left(\mathbb{E}\|X\|^{p} + 1 \right) < \infty$}
\end{array}
\]
proving (3.16). ~$\Box$

\vskip 0.3cm

The following recent result of Li, Qi, and Rosalsky [11] is used in the proof of
Theorem 2.3. It was proved by applying the new versions of the classical L\'{e}vy
and classical Hoffmann-J{\o}rgensen [6] inequalities established by Li and Rosalsky [12].

\vskip 0.3cm

\begin{theorem}
{\rm (Li, Qi, and Rosalsky [11])}. Let $q > 0$ and let $\{a_{n};~ n \geq 1 \}$
be a sequence of nonnegative real numbers such that
$\sum_{n=1}^{\infty} a_{n} < \infty$. Let $\{V_{k}; ~k \geq 1 \}$
be a sequence of independent symmetric $\mathbf{B}$-valued random variables.
Write
\[
b_{n} = \sum_{k=n}^{\infty} a_{k}, ~~n \geq 1
\]
and
\[
\alpha = \left \{
\begin{array}{ll}
2^{1-q}, & \mbox{if}~~0 < q \leq 1 \\
&\\
1, & \mbox{if}~~q > 1.
\end{array}
\right.
~~\mbox{and}~~
\beta = \left \{
\begin{array}{ll}
1, & \mbox{if}~~0 < q \leq 1 \\
&\\
2^{q - 1}, & \mbox{if}~~q > 1.
\end{array}
\right.
\]
Then, for all nonnegative real numbers $s$, $t$, and $u$, we have that
\[
\mathbb{P}\left( \sup_{n \geq 1} b_{n}
\left\|V_{n} \right\|^{q} > t \right)
\leq 2 \mathbb{P}\left( \sum_{n=1}^{\infty} a_{n} \left\|\sum_{i=1}^{n} V_{i} \right\|^{q}
> \frac{t}{\alpha} \right)
\]
and
\[
\begin{array}{ll}
& \mbox{$\displaystyle
\mathbb{P}\left( \sum_{n=1}^{\infty} a_{k}
\left\|\sum_{i=1}^{n} V_{i} \right\|^{q} > s + t + u \right)$}\\
&\\
& \mbox{$\displaystyle
\leq \mathbb{P}\left( \sup_{n \geq 1} b_{n}
\left\|V_{n} \right\|^{q} > \frac{s}{\beta^{2}} \right)
+ 4 \mathbb{P}\left( \sum_{n=1}^{\infty} a_{n} \left\|\sum_{i=1}^{n} V_{i} \right\|^{q}
> \frac{u}{\alpha \beta} \right) \mathbb{P}\left( \sum_{n=1}^{\infty}
a_{n} \left\|\sum_{i=1}^{n} V_{i} \right\|^{q} > \frac{t}{\alpha \beta^{2}} \right)$.}
\end{array}
\]
Furthermore, we have that
\[
\mathbb{E}\left(\sup_{n \geq 1} b_{n} \left\|V_{n} \right\|^{q} \right)
\leq 2 \alpha \mathbb{E}\left(\sum_{n=1}^{\infty} a_{k}
\left\|\sum_{i=1}^{n} V_{i} \right\|^{q} \right)
\]
and
\[
\mathbb{E}\left(\sum_{n=1}^{\infty} a_{k} \left\|\sum_{i=1}^{n} V_{i} \right\|^{q} \right)
\leq 6 (\alpha + \beta)^{3} \mathbb{E}\left(\sup_{n \geq 1} b_{n} \left\|V_{n} \right\|^{q} \right)
+ 6(\alpha + \beta)^{3} t_{0},
\]
where
\[
t_{0} = \inf \left\{t > 0;~\mathbb{P}\left( \sum_{n=1}^{\infty}
a_{n} \left\|\sum_{i=1}^{n} V_{i} \right\|^{q} > t \right) \leq 24^{-1} (\alpha + \beta)^{-3} \right\}.
\]
\end{theorem}

\vskip 0.3cm

\begin{lemma}
{\rm (Li, Qi, and Rosalsky [11])}. Let $({\bf E}, \mathcal{G})$ be a measurable linear
space and $g:~{\bf E} \rightarrow [0, \infty]$ be a measurable even function such that
for all $\mathbf{x}, \mathbf{y} \in {\bf E}$,
\[
g(\mathbf{x} + \mathbf{y}) \leq \beta \left(g(\mathbf{x}) + g(\mathbf{y}) \right),
\]
where $1 \leq \beta < \infty$ is a constant, depending only on the
function $g$. If $\mathbf{V}$ is an ${\bf E}$-valued random variable and
$\hat{\mathbf{V}}$ is a symmetrized version of $\mathbf{V}$
(i.e., $\hat{\mathbf{V}} = \mathbf{V} - \mathbf{V}^{\prime}$ where $\mathbf{V}^{\prime}$
is an independent copy of $\mathbf{V}$), then for all $t \geq 0$, we have that
\[
\mathbb{P}(g(\mathbf{V}) \leq t) \mathbb{E}g(\mathbf{V})
\leq \beta \mathbb{E} g(\hat{\mathbf{V}}) + \beta t
\]
and
\[
\mathbb{E}g(\hat{\mathbf{V}}) \leq 2 \beta \mathbb{E}g(\mathbf{V}).
\]
Moreover, if
\[
g(\mathbf{V}) < \infty ~~\mbox{a.s.},
\]
then
\[
\mathbb{E}g(\mathbf{V}) < \infty ~~\mbox{if and only if}~~\mathbb{E}g(\hat{\mathbf{V}}) < \infty.
\]
\end{lemma}

\vskip 0.3cm

\begin{lemma}
Let $1 \leq p < 2$ and let $\{X_{n};~n \geq 1\}$ be a sequence of independent copies
of a $\mathbf{B}$-valued random variable $X$ with $\mathbb{E}X = 0$ and
$\mathbb{E}\|X\|^{p} < \infty$. If $\mathbf{B}$ is a Banach space of Rademacher type $p$,
then
\begin{equation}
\sum_{n=1}^{\infty} \frac{1}{n}
\mathbb{E} \left(\frac{\left\|S^{(1)}_{n} - \mathbb{E}S^{(1)}_{n} \right\|^{p}}{n} \right) < \infty
~~\mbox{if and only if}~~\sum_{n=1}^{\infty} \frac{1}{n}
\mathbb{E} \left(\frac{\left\|U_{n} - \mathbb{E}U_{n} \right\|^{p}}{n} \right) < \infty.
\end{equation}
\end{lemma}

{\it Proof}~~Note that
\[
\left(S^{(1)}_{n} - \mathbb{E}S^{(1)}_{n}\right) - \left(U_{n} - \mathbb{E}U_{n} \right)
= \sum_{k=1}^{n}
\left(X_{k}I\{k < \|X_{k}\|^{p} \leq n \} - \mathbb{E} XI\{k < \|X\|^{p} \leq n \} \right),
~~n \geq 1.
\]
Then since $\mathbf{B}$ is a Banach space of Rademacher type $p$, we have that
\[
\begin{array}{lll}
\mbox{$\displaystyle
\mathbb{E} \left\|\left(S^{(1)}_{n} - \mathbb{E}S^{(1)}_{n}\right)
- \left(U_{n} - \mathbb{E}U_{n} \right)\right\|^{p}$}
&\leq&
\mbox{$\displaystyle C_{11} \sum_{k=1}^{n} \mathbb{E}
\left\|XI\{k < \|X\|^{p} \leq n \} - \mathbb{E} XI\{k < \|X\|^{p} \leq n \} \right|^{p}$}\\
&&\\
&\leq&
\mbox{$\displaystyle C_{12} \sum_{k=1}^{n}
\mathbb{E} \|X\|^{p} I\{k < \|X\|^{p} \leq n \}, ~~n \geq 1.$}
\end{array}
\]
Let $Y = \|X\|^{p}$. Then it follows from $\mathbb{E}Y < \infty$ (since $\mathbb{E}\|X\|^{p} < \infty$)
and the first conclusion of Lemma 3.5 (i.e., (3.13)) that
\[
\sum_{n=1}^{\infty}\frac{1}{n}
\mathbb{E}\left(\frac{\left\|\left(S^{(1)}_{n} - \mathbb{E}S^{(1)}_{n}\right)
- \left(U_{n} - \mathbb{E}U_{n} \right)\right\|^{p}}{n}\right)
\leq C_{12} \sum_{n=1}^{\infty} \frac{1}{n^{2}} \sum_{k=1}^{n} \mathbb{E}YI\{k < Y \leq n\} < \infty,
\]
which yields (3.17). ~$\Box$

\vskip 0.3cm

{\it Proof of Theorem 2.1}  To prove Theorem 2.1, we make the following simple observation.
Let $0 < p < q \leq 1$.  Let $ \{X_{n};~n \geq 1 \}$ be a sequence
of independent copies of a $\mathbf{B}$-valued random variable $X$
with $\mathbb{E}\|X\|^{p} < \infty$. Set $p_{1} = p/q$, $Y = \|X\|^{q}, ~Y_{n} = \|X_{n}\|^{q},
~n \geq 1$. Then $0 < p_{1} < 1$ and $\mathbb{E} Y^{p_{1}} < \infty$, and
\[
\sum_{n=1}^{\infty} \frac{1}{n} \left(\frac{\|S_{n}\|}{n^{1/p}} \right)^{q} \leq
\sum_{n=1}^{\infty} \frac{\sum_{k=1}^{n} \|X_{k}\|^{q}}{n^{1 + q/p}}
= \sum_{k=1}^{\infty} \left(\sum_{n=k}^{\infty} \frac{1}{n^{1+q/p}} \right) \|X_{k}\|^{q}
\leq \sum_{k=1}^{\infty} k^{-1/p_{1}} Y_{k} < \infty ~~\mbox{a.s.}
\]
(see Theorem 5.1.3 in Chow and Teicher [2, p. 118]). Theorem 2.1 follows immediately from
this observation together with Theorem 1.4 and Remark 1.2. ~$\Box$

\vskip 0.3cm

{\it Proof of Theorem 2.2} ({\bf Sufficiency}) Firstly we consider the case where $1 \leq q < p < 2$.
Since $\mathbb{E}X = 0$, we see that
\[
S_{n} = \left(U^{(1)}_{n} - \mathbb{E}U^{(1)}_{n} \right)
+ \left((S_{n} - U^{(1)}_{n}) - \mathbb{E}(S_{n} - U^{(1)}_{n}) \right), ~~n \geq 1
\]
so that, by Lemmas 3.3 and 3.4, (2.2) ensures (1.1) which implies (2.1).

Secondly we consider the case where $1 < p < q$. Since $\mathbf{B}$ is of stable type $p$,
the Maurey-Pisier [16] theorem asserts that it is also of stable type $p + \delta$ for
some $0 < \delta < q - p$. By Remark 1.2, (2.1) holds if we can show that
\begin{equation}
X \in SLLN(p, p + \delta); ~~\mbox{i.e.,}~
\sum_{n=1}^{\infty} \frac{1}{n} \left(\frac{\|S_{n}\|}{n^{1/p}} \right)^{p+\delta} < \infty
~~\mbox{a.s.}
\end{equation}
Since $\mathbb{E}X = 0$, we have that
\begin{equation}
\begin{array}{lll}
\mbox{$\displaystyle S_{n}$}
&=&
\mbox{$\displaystyle
\sum_{k=1}^{n}X_{k}I\{\|X_{k}\|^{p} \leq n\} + \sum_{k=1}^{n}X_{k}I\{\|X_{k}\|^{p} > n\}$}\\
&&\\
&=&
\mbox{$\displaystyle
\left(U_{n} - \mathbb{E}U_{n} \right)
- n \mathbb{E}XI\{\|X\|^{p} > n\} + \sum_{k=1}^{n}X_{k}I\{\|X_{k}\|^{p} > n\}, ~~n \geq 1.$}
\end{array}
\end{equation}
It is easy to see that
\[
\left\{ \max_{1 \leq k \leq n} \|X_{k}\|^{p} > n ~~\mbox{i.o.} (n) \right\}
= \left\{ \|X_{n}\|^{p} > n ~~\mbox{i.o.} (n) \right\}.
\]
Since $\{X_{n};~n \geq 1 \}$ is a sequence of independent copies of $\mathbf{B}$-valued
random variable $X$ with $\mathbb{E}\|X\|^{p} < \infty$, it follows from the Borel-Cantelli
lemma that
\[
\mathbb{P}\left( \|X_{n}\|^{p} > n ~\mbox{i.o.} (n) \right) = 0
\]
and hence
\begin{equation}
\mathbb{P}\left(\max_{1 \leq k \leq n} \|X_{k}\|^{p} > n ~\mbox{i.o.} (n) \right)
= 0,
\end{equation}
which ensures that
\begin{equation}
\sum_{n=1}^{\infty} \frac{1}{n}
\left(\frac{\left\|\sum_{k=1}^{n}X_{k}I\{\|X_{k}\|^{p} > n\} \right\|}{n^{1/p}} \right)^{p+\delta}
< \infty~~\mbox{a.s.}
\end{equation}
Note that $1 < p < 2$ and $\mathbb{E}\|X\|^{p} < \infty$ imply that
\[
\frac{n\mathbb{E}\|X\|I\{\|X\|^{p} > n\}}{n^{1/p}}
 \leq \frac{n^{1/p}\mathbb{E}\|X\|^{p}}{n^{1/p}} = \mathbb{E}\|X\|^{p}, ~~n \geq 1.
\]
Thus, by (3.16) of Lemma 3.5, we have that
\begin{equation}
\begin{array}{ll}
& \mbox{$\displaystyle
\sum_{n=1}^{\infty} \frac{1}{n}
\left(\frac{\|n \mathbb{E}XI\{\|X\|^{p} > n\}\|}{n^{1/p}} \right)^{p+\delta}$}\\
&\\
&\mbox{$\displaystyle
= \sum_{n=1}^{\infty}
\frac{1}{n} \left(\frac{\|n \mathbb{E}XI\{\|X\|^{p} > n\}\|}{n^{1/p}} \right)^{p}
\left(\frac{\|n \mathbb{E}XI\{\|X\|^{p} > n\}\|}{n^{1/p}} \right)^{\delta}$}\\
&\\
& \mbox{$\displaystyle \leq \left(\mathbb{E}\|X\|^{p} \right)^{\delta}
\sum_{n=1}^{\infty} \frac{\left(\mathbb{E}\|X\| I\{\|X\|^{p} > n\} \right)^{p}}{n^{2-p}}$}\\
&\\
&\mbox{$\displaystyle < \infty.$}
\end{array}
\end{equation}
Since $\mathbf{B}$ is also of Rademacher type $p+\delta$, we get that
\[
\begin{array}{lll}
\mbox{$\displaystyle
\mathbb{E}\left\|U_{n} - \mathbb{E}U_{n} \right\|^{p + \delta}$}
&\leq&
\mbox{$\displaystyle
C_{13} n \mathbb{E}
\left\|X I\{\|X\|^{p} > n\} - \mathbb{E} X I\{\|X\|^{p} > n\} \right\|^{p+\delta}$}\\
&&\\
&\leq&
\mbox{$\displaystyle C_{14} \mathbb{E}\|X\|^{p+\delta} I\{\|X\|^{p} > n\}, ~~n \geq 1.$}
\end{array}
\]
Thus, by (3.15) of Lemma 3.5, we have that
\[
\sum_{n=1}^{\infty} \frac{1}{n}
\mathbb{E}\left(\frac{\left\|U_{n} - \mathbb{E}U_{n} \right\|}{n^{1/p}} \right)^{p+\delta}
\leq C_{14} \sum_{n=1}^{\infty}
\frac{\mathbb{E}\|X\|^{p+\delta} I\{\|X\|^{p} > n\}}{n^{1 + \delta/p}} < \infty
\]
and hence
\[
\sum_{n=1}^{\infty} \frac{1}{n}
\left(\frac{\left\|U_{n} - \mathbb{E}U_{n} \right\|}{n^{1/p}} \right)^{p+\delta}
< \infty ~~\mbox{a.s.,}
\]
which, together with (3.19), (3.21), and (3.22), ensures (3.18).

Lastly we consider the case where $1 < q = p < 2$. Since $\mathbb{E}\|X\|^{p} < \infty$,
we have that
\[
\lim_{n \rightarrow \infty} \frac{u_{n}^{p}}{n} = 0;
~~\mbox{i.e.,}~~\lim_{n \rightarrow \infty} \frac{u_{n}}{n^{1/p}} = 0.
\]
Hence we can assume, without loss of generality, that $u_{n} < n^{1/p}$ for all $n \geq 1$.
Since $\mathbb{E}X = 0$, we have that, for $n \geq 1$,
\begin{equation}
\begin{array}{lll}
\mbox{$\displaystyle S_{n} $}
&=&
\mbox{$\displaystyle
\sum_{k=1}^{n} X_{k}I\{\|X_{k}\| \leq u_{n} \} + \sum_{k=1}^{n} X_{k}I\{u_{n} < \|X_{k}\| \leq n^{1/p} \}
+ \sum_{k=1}^{n} X_{k}I\{\|X_{k}\| > n^{1/p} \}$}\\
&&\\
&=& \mbox{$\displaystyle
\left(U^{(1)}_{n} - \mathbb{E}U^{(1)}_{n} \right) + \left(U^{(2)}_{n} - \mathbb{E}U^{(2)}_{n} \right)
- n \mathbb{E}X I\left\{\|X\|^{p} > n\right\}
+ \sum_{k=1}^{n} X_{k}I\left\{\|X_{k}\|^{p} > n \right\}.$}\\
\end{array}
\end{equation}
Since (3.20) follows from $\mathbb{E}\|X\|^{p} < \infty$, we see that
\begin{equation}
\sum_{n=1}^{\infty} \frac{1}{n}
\left(\frac{\left\|\sum_{k=1}^{n}X_{k}I\{\|X_{k}\|^{p} > n\} \right\|}{n^{1/p}} \right)^{p}
< \infty~~\mbox{a.s.}
\end{equation}
Since $p > 1$ and $\mathbb{E}\|X\|^{p} < \infty$, it follows from (3.16) of Lemma 3.5 that
\begin{equation}
\sum_{n=1}^{\infty} \frac{1}{n}
\left(\frac{\left\|n \mathbb{E}X I\left\{\|X\|^{p} > n\right\} \right\|}{n^{1/p}} \right)^{p}
\leq \sum_{n=1}^{\infty}
\frac{\left(\mathbb{E}\|X\| I\left\{\|X\|^{p} > n\right\} \right)^{p}}{n^{2-p}} < \infty.
\end{equation}
Since $\mathbb{E}\|X\|^{p} < \infty$ and $\mathbf{B}$ is a Banach space of stable type
$p \in (1, 2)$, by Remark 3.2, (3.12) holds, which ensures that
\begin{equation}
\sum_{n=1}^{\infty} \frac{1}{n}
\left(\frac{\left\|U^{(1)}_{n} - \mathbb{E}U^{(1)}_{n} \right\|}{n^{1/p}} \right)^{p}
< \infty ~~\mbox{a.s.}
\end{equation}
Since $\mathbf{B}$ is also a Banach space of Rademacher type $p$, we have that,
for all $n \geq 1$,
\[
\begin{array}{lll}
\mbox{$\displaystyle \mathbb{E}\|U_{n}^{(2)} -
\mathbb{E}U_{n}^{(2)}\|^{p}$} &\leq& \mbox{$\displaystyle
C_{15}\sum_{k=1}^{n}\mathbb{E}\left\|X_{k}I\left\{u_{n} < \|X_{k}\| \leq n^{1/p}
\right\} - \mathbb{E}\left(XI\left\{u_{n} < \|X\| \leq n^{1/p}
\right\} \right) \right\|^{p}$}\\
&&\\
&\leq& \mbox{$\displaystyle 2C_{15}
n \mathbb{E}\left(\|X\|^{p}I\left\{u_{n}^{p} < \|X\|^{p} \leq n \right\} \right)$}\\
&&\\
&=& \mbox{$\displaystyle 2C_{15} n\int_{u_{n}^{p}}^{n} t
d\mathbb{P}\left(\|X\|^{p} \leq t \right)$}\\
&&\\
&\leq& \mbox{$\displaystyle 2C_{15} nu_{n}^{p}
\mathbb{P}\left(\|X\|^{p} > u_{n}^{p}\right)
+ 2C_{15} n\int_{u_{n}^{p}}^{n} \mathbb{P}\left(\|X\|^{p} > t \right) dt$}\\
&&\\
&\leq& \mbox{$\displaystyle 2C_{15}u_{n}^{p} + 2C_{15}n\int_{u_{n}^{p}}^{n}
\mathbb{P}\left(\|X\|^{p} > t \right) dt.$}
\end{array}
\]
Now (3.14) holds by Lemma 3.5. Thus it follows from (3.14) and
(2.2) that
\[
\sum_{n=1}^{\infty} \frac{\mathbb{E}\|U_{n}^{(2)}
- \mathbb{E}U_{n}^{(2)}\|^{p}}{n^{2}} < \infty,
\]
which ensures that
\begin{equation}
\sum_{n=1}^{\infty} \frac{\|U_{n}^{(2)}
- \mathbb{E}U_{n}^{(2)}\|^{p}}{n^{2}} < \infty ~~\mbox{a.s.}
\end{equation}
Combining (3.23)-(3.27), we conclude that (2.1) holds for $q = p$.
The proof of the sufficiency half of Theorem 2.2 is complete. ~$\Box$

\vskip 0.3cm

{\it Proof of Theorem 2.2} ~({\bf Necessity}) For the case where $q \neq p$, by Theorem 1.4,
we see that (2.2) follows immediately from (2.1).

We now consider the case where $q = p$. By Theorem 1.4, (2.1) implies that $\mathbb{E}X = 0$
and $\mathbb{E}\|X\|^{p} < \infty$. Hence we can assume, without loss of generality, that
$u_{n}^{p} < n$ for all $n \geq 1$. We thus only need to show that (2.1) (with $q = p$)
implies that
\begin{equation}
\sum_{n=1}^{\infty} \frac{\int_{u_{n}^{p}}^{n} \mathbb{P} \left(\|X\|^{p} > t\right)dt}{n}
< \infty.
\end{equation}
To see this, let $\{X^{\prime}, ~X^{\prime}_{n};~n \geq 1 \}$ be an independent copy
of $\{X, ~X_{n};~n \geq 1 \}$. Let
\[
V_{n} = \left(X_{n}I\left\{\|X_{n}\|^{p} \leq n \right\}
- X^{\prime}_{n}I\left\{\|X^{\prime}_{n}\|^{p} \leq n \right\} \right)~~\mbox{and}~~
\hat{S}^{(1)}_{n} = \sum_{k=1}^{n}V_{k}, ~~n \geq 1.
\]
Then $\{V_{n};~n \geq 1 \}$ is a sequence of independent symmetric $\mathbf{B}$-valued
random variables. By the Borel-Cantelli lemma, it follows from $\mathbb{E}\|X\|^{p} < \infty$
that
\[
\mathbb{P} \left(\|X_{n}\|^{p} > n ~\mbox{i.o.} (n) \right) = 0,
\]
which ensures that
\[
S^{(2)}_{n} = \sum_{k=1}^{n} X_{k}I\left\{\|X_{k}\|^{p} > k \right\}
= {\it O}(1) ~~\mbox{a.s. as}~n \rightarrow \infty
\]
and hence
\[
\sum_{n=1}^{\infty} \frac{1}{n} \left(\frac{\|S^{(2)}_{n} \|}{n^{1/p}} \right)^{p}
= \sum_{n=1}^{\infty} \frac{\|S^{(2)}_{n} \|^{p}}{n^{2}}
< \infty ~~\mbox{a.s.}
\]
Note that
\[
\|S_{n}^{(1)}\| \leq \|S_{n}\| + \|S_{n}^{(2)} \|, ~~n \geq 1.
\]
It thus follows from (2.1) (with $q = p$) that
\begin{equation}
\sum_{n=1}^{\infty} \frac{\|S^{(1)}_{n} \|^{p}}{n^{2}}
< \infty~~\mbox{a.s.}
\end{equation}
and hence
\begin{equation}
\sum_{n=1}^{\infty} \frac{\|\hat{S}^{(1)}_{n} \|^{p}}{n^{2}}
< \infty~~\mbox{a.s.}
\end{equation}
Let $a_{n} = 1/n^{2}, ~n \geq 1$. Then
\[
b_{n} = \sum_{k=n}^{\infty} a_{k} = \sum_{k=n}^{\infty}
\frac{1}{n^{2}} \leq \frac{2}{n}, ~~n \geq 1
\]
and hence
\[
\sup_{n \geq 1} b_{n} \|V_{n}\|^{p}
\leq \sup_{n \geq 1} \frac{2}{n} \left(2n^{1/p} \right)^{p} = 2^{p+1} ~~\mbox{a.s.}
\]
We thus have that
\begin{equation}
\mathbb{E} \left(\sup_{n \geq 1} b_{n} \|V_{n}\|^{p} \right) < \infty.
\end{equation}
By Theorem 3.1, we conclude from (3.30) and (3.31) that
\[
\sum_{n=1}^{\infty} \frac{\mathbb{E}\|\hat{S}^{(1)}_{n} \|^{p}}{n^{2}}
< \infty;
\]
that is,
\begin{equation}
\mathbb{E} \left(\sum_{n=1}^{\infty} \frac{\|\hat{S}^{(1)}_{n} \|^{p}}{n^{2}} \right)
< \infty.
\end{equation}
By Lemma 3.6, it follows from (3.29) and (3.32) that
\[
\mathbb{E} \left(\sum_{n=1}^{\infty} \frac{\|S^{(1)}_{n} \|^{p}}{n^{2}} \right)
< \infty;
\]
that is,
\begin{equation}
\sum_{n=1}^{\infty} \frac{\mathbb{E}\|S^{(1)}_{n} \|^{p}}{n^{2}}
< \infty.
\end{equation}
Since $1 < p < 2$, applying (2.5) of Ledoux and Talagrand [9, p. 46], (3.33)
ensures that
\[
\sum_{n=1}^{\infty} \frac{\|\mathbb{E}S^{(1)}_{n} \|^{p}}{n^{2}}
< \infty
\]
which, together with (3.33), gives
\[
\sum_{n=1}^{\infty} \frac{\mathbb{E}\|S^{(1)}_{n} - \mathbb{E}S^{(1)}_{n}\|^{p}}{n^{2}}
< \infty.
\]
By Lemma 3.7, this is equivalent to
\begin{equation}
\sum_{n=1}^{\infty} \frac{\mathbb{E}\|U_{n} - \mathbb{E}U_{n}\|^{p}}{n^{2}}
< \infty.
\end{equation}
Since $\|U_{n}^{(2)} - \mathbb{E}U_{n}^{(2)}\| \leq \|U_{n} -
\mathbb{E}U_{n}\| + \|U_{n}^{(1)} - \mathbb{E}U_{n}^{(1)}\|, ~n \geq
1$ and {\bf B} is of stable type $p$ where $1 < p < 2$, it follows from
Remark 3.2 and (3.34) that
\begin{equation}
\sum_{n=1}^{\infty} \frac{\mathbb{E}\|U_{n}^{(2)}
- \mathbb{E}U_{n}^{(2)}\|^{p}}{n^{2}} < \infty.
\end{equation}
By Lemma 3.1 {\bf (ii)} of Li, Qi, and Rosalsky [10],
\[
\begin{array}{ll}
& \mbox{$\displaystyle
\mathbb{E} \max_{1 \leq k \leq n}
\left\| X_{k}I\{u_{n}^{p} < \|X_{k}\|^{p} \leq n\}
- \mathbb{E}XI\{u_{n}^{p} < \|X\|^{p} \leq n\} \right\|^{p}$}\\
&\\
& \mbox{$\displaystyle
= \mathbb{E} \left(\max_{1 \leq k \leq n} \left\| X_{k}I\{u_{n}^{p} < \|X_{k}\|^{p}
\leq n\} - \mathbb{E}XI\{u_{n}^{p} < \|X\|^{p} \leq n\} \right\|\right)^{p}$}\\
&\\
& \mbox{$\displaystyle \leq 2^{p+1}
\mathbb{E}\|U_{n}^{(2)} - \mathbb{E}U_{n}^{(2)} \|^{p}, ~n \geq 1.$}
\end{array}
\]
It thus follows from (3.35) and $\mathbb{E}\|X\|^{p} < \infty$ that
\[
\begin{array}{ll}
& \mbox{$\displaystyle \sum_{n=1}^{\infty}
\frac{\mathbb{E} \max_{1 \leq k \leq n}
\left\|X_{k}I\{u_{n}^{p} < \|X_{k}\|^{p} \leq n\}\right\|^{p}}{n^{2}}$}\\
&\\
& \mbox{$\displaystyle \leq 2^{2p} \sum_{n=1}^{\infty}
\frac{\mathbb{E}\|U_{n}^{(2)} - \mathbb{E}U_{n}^{(2)} \|^{p}}{n^{2}}
+ 2^{p-1}\sum_{n=1}^{\infty}\frac{\mathbb{E}\|X\|^{p}}{n^{2}}$}\\
&\\
& \mbox{$\displaystyle < \infty,$}
\end{array}
\]
and hence, by Lemma 5.4 of Li, Qi, and Rosalsky [10], noting that
$\mathbb{P}\left(\|X\|^{p} > u_{n}^{p} \right) \leq n^{-1}, ~n \geq 1$,
we get that
\begin{equation}
\sum_{n=1}^{\infty}
\frac{\mathbb{E}\|X\|^{p}I\left\{u_{n}^{p} < \|X\|^{p} \leq n \right\}}{n}
< \infty.
\end{equation}
Using partial integration, one can easily see that
\begin{equation}
\left|\mathbb{E}\|X\|^{p}I\left\{u_{n}^{p} < \|X\|^{p} \leq n \right\}
- \int_{u_{n}^{p}}^{n} \mathbb{P}\left(\|X\|^{p} > t \right)dt \right|
\leq \frac{u_{n}^{p}}{n}
+ n\mathbb{P}\left(\|X\|^{p} > n \right), ~~ n \geq 1.
\end{equation}
Since $\mathbb{E}\|X\|^{p} < \infty$, we have
\begin{equation}
\sum_{n=1}^{\infty} \frac{n\mathbb{P}\left(\|X\|^{p} > n \right)}{n} =
\sum_{n=1}^{n} \mathbb{P}\left(\|X\|^{p} > n \right) < \infty,
\end{equation}
and, by Lemma 3.5, (3.14) holds. We thus see that (3.28) follows from
(3.36), (3.37), (3.38), and (3.14) thereby completing the proof of
the necessity half of Theorem 2.2.~ $\Box$

\vskip 0.3cm

{\it Proof of Theorem 2.3}~~We only need to consider the case where $q > 1$ since
for the case where $q = 1$, Theorem 2.3 is Theorem 2.3 of Li, Qi, and Rosalsky [10].
Note that
\[
\left\{\max_{1 \leq k \leq n}\|X_{k}\| > n ~\mbox{i.o.} (n) \right\}
= \left\{\|X_{n}\| > n ~\mbox{i.o.} (n) \right\}
\]
and for $p =1$,
\[
U_{n} = \sum_{k=1}^{n} X_{k}I\left\{\|X_{k}\| \leq n \right\}, ~~n \geq 1.
\]
By the Borel-Cantelli lemma, it thus follows from $\mathbb{E}\|X\| < \infty$ that
\[
\mathbb{P} \left(\max_{1 \leq k \leq n}\|X_{k}\| > n ~\mbox{i.o.} (n) \right) = 0
\]
and hence
\[
\mathbb{P} \left(S_{n} - U_{n} \neq 0 ~\mbox{i.o.} (n) \right) = 0,
\]
which ensures that
\begin{equation}
\sum_{n=1}^{\infty} \frac{1}{n}
\left( \frac{\left\|S_{n} - U_{n} \right\|}{n} \right)^{q}
< \infty ~~\mbox{a.s.}
\end{equation}
and by the Mourier [17] SLLN, it follows from (2.4) that
\begin{equation}
\lim_{n \rightarrow \infty} \frac{U_{n}- n \mathbb{E}(XI\{\|X\| \leq n\})}{n}
= \lim_{n \rightarrow \infty} \left(\frac{S_{n}}{n} - \mathbb{E}(XI\{\|X\| \leq n\}) \right)
- \lim_{n \rightarrow \infty} \frac{S_{n} - U_{n}}{n} = 0 ~~\mbox{a.s.}
\end{equation}
We now show that
\begin{equation}
\sum_{n=1}^{\infty} \frac{1}{n}
\left( \frac{\left\|U_{n} - n \mathbb{E}(XI\{\|X\| \leq n\}) \right\|}{n} \right)^{q}
< \infty ~~\mbox{a.s.}
\end{equation}
Since $\mathbf{B}$ is of stable type $1$, the Maurey-Pisier [16]
theorem asserts that it is also of stable type $1+\delta$ for some $0 < \delta < q-1$
and hence
\[
\mathbb{E} \left\|U_{n} - n \mathbb{E}(XI\{\|X\| \leq n\}) \right\|^{1 + \delta}
\leq C_{16} \mathbb{E}\left(\|X\|^{1+\delta}I\{\|X\| \leq n\} \right), ~~n \geq 1.
\]
Thus, by (3.15) (with $p = 1$) of Lemma 3.5, we conclude that
\[
\sum_{n=1}^{\infty} \frac{1}{n} \mathbb{E}
\left( \frac{\left\|U_{n} - n \mathbb{E}(XI\{\|X\| \leq n\}) \right\|}{n} \right)^{1+ \delta}
< \infty
\]
and hence
\[
\sum_{n=1}^{\infty} \frac{1}{n}
\left( \frac{\left\|U_{n} - n \mathbb{E}(XI\{\|X\| \leq n\}) \right\|}{n} \right)^{1+ \delta}
< \infty ~~\mbox{a.s.},
\]
which, together with (3.40), ensures that (3.41) holds since $q > 1 + \delta$. Note that
\[
S_{n} = \left(S_{n} - U_{n} \right) + \left(U_{n} - n \mathbb{E}(XI\{\|X\| \leq n\}) \right)
+ n\mathbb{E}(XI\{\|X\| \leq n\}), ~~n \geq 1.
\]
We thus see that (2.3) (with $q > 1$) follows from (3.39), (3.41), and the second half of (2.4)
(with $q > 1$).

Conversely, by Theorem 1.4 and the Mourier [17] SLLN, it follows from (2.3) that $\mathbb{E}X = 0$
and $\mathbb{E}\|X\| < \infty$ and hence (3.39) and (3.41) (since $\mathbf{B}$ is of stable type $1$)
hold. Note that
\[
n \mathbb{E}(XI\{\|X\| \leq n\}) = S_{n}
- \left(S_{n} - U_{n} \right) - \left(U_{n} - n \mathbb{E}(XI\{\|X\| \leq n\})\right), ~~n \geq 1.
\]
It thus follows from (2.3), (3.39), and (3.41) that
\[
\sum_{n=1}^{\infty} \frac{\|\mathbb{E}(XI\{\|X\| \leq n\})\|^{q}}{n}
= \sum_{n=1}^{\infty} \frac{1}{n} \left(\frac{\|n \mathbb{E}(XI\{\|X\| \leq n\})\|}{n}\right)^{q}
< \infty
\]
and hence (2.4) holds (with $q > 1$). The proof of Theorem 2.3 is complete. ~$\Box$

\section{Three Examples}

Li, Qi, and Rosalsky [10] provided three examples (see, Examples 5.1, 5.2, and 5.3 of Li, Qi, and
Rosalsky [10]) for illustrating the necessary and sufficient conditions that they obtained for
(2.3) for the case where $q = 1$. In this section we provide three examples to illustrate our
Theorems 1.4, 2.2, and 2.3.

\vskip 0.3cm

\begin{example}
Let $1 < r < p < 2$ and let $X$ be a real-valued symmetric random variable such that
\[
\mathbb{P}(X = 0) = b ~~\mbox{and}~~\mathbb{P}(|X| > t)
= \int_{t}^{\infty} \frac{1}{x^{p+1} \ln^{r}t}dt, ~~t \geq e,
\]
where $b = 1 - \int_{e}^{\infty} \frac{1}{x^{p+1} \ln^{r}x}dx$. Then
\[
\mathbb{P}(|X| > t) \sim \frac{1}{p t^{p}\ln^{r}t} ~~\mbox{as}~~t \rightarrow \infty
\]
and hence, for $1 \leq q < p$,
\[
\mathbb{P}^{q/p}\left(|X|^{q} > t \right) = \mathbb{P}^{q/p}\left(|X| > t^{1/q}\right) \sim
(q^{r}/p)^{q/p} t^{-1} (\ln t)^{-rq/p}~~\mbox{as}~~t \rightarrow \infty.
\]
We then see that
\[
\int_{0}^{\infty} \mathbb{P}^{q/p}\left(|X|^{q} > t \right)dt
\left\{
\begin{array}{ll}
< \infty & \mbox{if}~p/r < q < p,\\
&\\
= \infty & \mbox{if}~1 \leq q \leq p/r.
\end{array}
\right.
\]
It is also easy to check that
\[
\mathbb{E}|X|^{p} \ln(1 + |X|) = \infty
~~\mbox{and}~~\mathbb{E}|X|^{q} = \infty ~~\mbox{for all}~ q > p.
\]
By Theorem 2.2 and Remark 1.2, for this example, $X \in SLLN(p, q)$ if and only if $p/r < q < \infty$.
However, by Corollary 2.2, (1.1) holds if and only if $p/r < q < p$. This means that, if (1.1)
holds for some $q = q_{1} > 0$, one cannot conclude that (1.1) holds for either $0 < q < q_{1}$
or $q > q_{1}$.
\end{example}

\vskip 0.3cm

\begin{example}
Let $1 < p < 2$ and let $X$ be a real-valued symmetric random variable with density
function
\[
f(x) = \frac{b}{|x|^{p+1}(\ln|x|) (\ln\ln|x|)^{2}}I\{|x| > 3\},
\]
where $0 < b < \infty$ is such that $\int_{-\infty}^{\infty}f(x)dx =
1$. Clearly, we have that
\[
\mathbb{E}X = 0 ~~\mbox{and}~~\mathbb{E}|X|^{p} < \infty.
\]
Since
\[
\mathbb{P}(|X| > x) \sim \frac{2b/p}{x^{p} (\ln x) (\ln \ln
x)^{2}}~~\mbox{as}~ x \rightarrow \infty,
\]
we see that
\[
u_{n} \sim \frac{(2bn)^{1/p}}{(\ln n)^{1/p} (\ln \ln n)^{2/p}} ~~\mbox{as}~ n
\rightarrow \infty
\]
and hence, for all sufficiently large $n$,
\[
\begin{array}{lll}
\mbox{$\displaystyle \int_{u_{n}^{p}}^{n} \mathbb{P}\left(|X|^{p} > t \right) dt$}
&=& \mbox{$\displaystyle \int_{u_{n}^{p}}^{n} \mathbb{P}\left(|X| > t^{1/p} \right) dt$}\\
&&\\
&\geq& \mbox{$\displaystyle \int_{\frac{3bn}{(\ln n) (\ln \ln
n)^{2}}}^{n} \frac{b}{t (\ln t) (\ln \ln t)^{2}}dt$}\\
&&\\
&\geq& \mbox{$\displaystyle \frac{b}{(\ln n)(\ln \ln n)^{2}}
\int_{\frac{3bn}{(\ln n) (\ln \ln
n)^{2}}}^{n} \frac{1}{t} dt$}\\
&&\\
&\sim& \mbox{$\displaystyle \frac{b}{(\ln n) (\ln \ln
n)}~~\mbox{as}~ n \rightarrow \infty.$}
\end{array}
\]
Note that
\[
\sum_{n=3}^{\infty}\frac{b}{n(\ln n)(\ln \ln n)} = \infty
\]
and so
\[
\sum_{n=1}^{\infty} \frac{\int_{\min\{u_{n}^{p},
n\}}^{n} \mathbb{P}\left(|X|^{p} > t \right) dt}{n} = \infty.
\]
By Theorem 2.2 and Remark 1.2, we thus conclude that $X \notin SLLN(p, q)$
for this example for all $0 < q \leq p$.
\end{example}

\vskip 0.3cm

Let $1 < p < 2$ and let $\{X_{n};~n \geq 1\}$ be a sequence of independent
copies of a symmetric real-valued random variable $X$. Then, by either
Theorem 2.2 or Theorem 2.4, the following three statements are equivalent:
\begin{description}
\item
\quad {\bf (i)} ~~~~$\displaystyle \mathbb{E}X = 0
~~\mbox{and}~~\mathbb{E}|X|^{p} < \infty;$

\item
\quad {\bf (ii)} ~~~$X \in SLLN(p, q)$ for some $q > p$;

\item
\quad {\bf (iii)} ~~$X \in SLLN(p, q)$ for all $q > p$.
\end{description}
However, the following example says that this is not true when $p = 1$.

\vskip 0.3cm

\begin{example}
Let $X$ be a real-valued random variable such that
\[
\mathbb{P}\left(X = - \frac{1}{1-a} \right) = 1 -a~~ \mbox{and}~~
\mathbb{P}(X > x) = \int_{x}^{\infty} \frac{1}{t^{2}(\ln t) (\ln \ln t)^{2}} dt,
~~x \geq e^{e}
\]
where $a = \int_{e}^{\infty} \frac{1}{t^{2}(\ln t) (\ln \ln t)^{2}} dt$. Then
\[
\mathbb{E}X = 0, ~~\mathbb{E}|X| < \infty,
\]
and, for all sufficiently large $n$,
\[
\mathbb{E}XI\{|X| \leq  n\} = - \mathbb{E}XI\{|X| > n\} = -
\int_{n}^{\infty} \frac{1}{t (\ln t) (\ln \ln t)^{2}}dt = -\frac{1}{\ln \ln n}.
\]
Note that
\[
\sum_{n=2}^{\infty}\frac{1}{n (\ln \ln n)^{q}} = \infty ~~\mbox{for all} ~ q > 1.
\]
Thus for this example, by either Theorem 2.3 or Theorem 2.4, $X \notin SLLN(1, q)$
for all $q > 1$.
\end{example}

\vskip 0.5cm

\noindent
{\bf Acknowledgments}\\

\noindent The research of Deli Li was partially supported
by a grant from the Natural Sciences and Engineering Research Council of
Canada and the research of Yongcheng Qi was partially supported by
NSF Grant DMS-1005345.

\vskip 0.5cm

{\bf References}

\begin{enumerate}

\item Azlarov, T. A., Volodin, N. A.: Laws of large numbers for
identically distributed Banach-space valued random variables. Teor.
Veroyatnost. i Primenen. {\bf 26}, ~584-590 (1981), in Russian.
English translation in Theory Probab. Appl. {\bf 26}, 573-580
(1981).

\item Chow, Y.S., Teicher, H.: Probability Theory: Independence,
Interchangeability, Martingales, 3rd ed. Springer-Verlag, New York
(1997).

\item de Acosta, A.: Inequalities for {\it B}-valued random vectors
with applications to the law of large numbers. Ann. Probab. {\bf 9},
157-161 (1981).

\item Hechner, F.: Comportement asymptotique de sommes de
Ces\`{a}ro al\'{e}atoires. C. R. Math. Acad. Sci. Paris {\bf 345}, 705-708
(2007).

\item Hechner, F., Heinkel, B.: The Marcinkiewicz-Zygmund LLN in
Banach spaces: A generalized martingale approach. J. Theor. Probab.
{\bf 23}, 509-522 (2010).

\item Hoffmann-J{\o}rgensen, J.: Sums of independent Banach
space valued random variables. Studia Math. {\bf 52}, 159-186
(1974).

\item Hoffmann-J{\o}rgensen, J., Pisier, G.: The law of large
numbers and the central limit theorem in Banach spaces. Ann. Probab.
{\bf 4}, 587-599 (1976).

\item Kolmogoroff, A.: Sur la loi forte des grands nombres. C. R.
Acad. Sci. Paris S\'{e}r. Math. {\bf 191}, 910-912 (1930).

\item Ledoux, M., Talagrand, M.:  Probability in Banach Spaces:
Isoperimetry and Processes. Springer-Verlag, Berlin (1991).

\item Li, D., Qi, Y., Rosalsky, A.: A refinement of the
Kolmogorov-Marcinkiewicz-Zygmund strong law of large numbers.
J. Theoret. Probab. {\bf 24}, 1130-1156 (2011).

\item Li, D., Qi, Y., Rosalsky, A.: An extension of a theorem of Hechner and
Heinkel. manuscript (2012).

\item Li, D., Rosalsky, A.: New versions of some classical stochastic
inequalities. Stochastic Anal. Appl. {\bf 31}, (2013) (to appear).

\item Marcinkiewicz, J., Zygmund, A.: Sur les fonctions
ind\'{e}pendantes. Fund. Math. {\bf 29}, 60-90 (1937).

\item Marcus, M. B., Pisier, G.: Characterizations of almost surely
continuous $p$-stable random Fourier series and strongly stationary
processes. Acta Math. {\bf 152}, 245-301 (1984).

\item Marcus, M. B., Woyczy\'{n}ski, W. A.: Stable measures and
central limit theorems in spaces of stable type. Trans. Amer. Math.
Soc. {\bf 251}, 71-102 (1979).

\item Maurey, B., Pisier, G.: S\'{e}ries de variables al\'{e}atoires
vectorielles ind\'{e}pendantes et propri\'{e}t\'{e}s
g\'{e}om\'{e}triques des espaces de Banach. Studia Math. {\bf 58},
45-90 (1976).

\item Mourier, E.: El\'{e}ments al\'{e}atoires dans un espace de
Banach. Ann. Inst. H. Poincar\'{e} {\bf 13}, 161-244 (1953).

\item Pisier, G.: Probabilistic methods in the geometry of Banach
spaces, in Probability and Analysis, Lectures given at the 1st 1985
Session of the Centro Internazionale Matematico Estivo (C.I.M.E.),
Lecture Notes in Mathematics, Vol. {\bf 1206}, 167-241,
Springer-Verlag, Berlin (1986).

\item Rosi\'{n}ski, J.: Remarks on Banach spaces of stable type.
Probab. Math. Statist. {\bf 1}, 67-71 (1980).

\item Woyczy\'{n}ski, W. A.: Geometry and martingales in Banach
spaces-Part II: Independent increments, in Probability on Banach
Spaces (Edited by J. Kuelbs), Advances in Probability and Related
Topics Vol. {\bf 4} (Edited by P. Ney), 267-517, Marcel Dekker, New
York (1978).

\end{enumerate}

\end{document}